\documentclass[11pt,reqno]{amsproc}
\usepackage[utf8]{inputenc}

\allowdisplaybreaks
\numberwithin{equation}{section}

\usepackage{color}
\usepackage{graphicx}
\usepackage{fullpage}
\usepackage{enumerate}
\usepackage{tcolorbox}

\usepackage[semicolon,square,authoryear,sort]{natbib}
\usepackage{amsmath}
\usepackage{xcolor}
\usepackage{epsfig}
\usepackage{subcaption}

\usepackage[debug=false, colorlinks=true, pdfstartview=FitV,
linkcolor=blue, citecolor=blue, urlcolor=blue]{hyperref}

\usepackage{tikz}
\usetikzlibrary{shapes.geometric, arrows}
\tikzstyle{startstop} = [rectangle, rounded corners, minimum width=3cm, minimum height=1cm,text centered, text width=5.5cm, draw=black, fill=red!10]
\tikzstyle{process} = [rectangle, rounded corners, minimum width=3cm, minimum height=1cm,text centered, text width = 4cm, draw=black, fill=orange!10]
\tikzstyle{arrow} = [thick,->,>=stealth]
\tikzstyle{innerWhite} = [semithick, white,line width=1.4pt, shorten >= 4.5pt]

\usepackage[doipre={DOI:~}]{uri}

\usepackage{makecell}
\usepackage[section]{placeins}

\usepackage{morefloats}
\usepackage{multirow}

\newlength{\drop}
\definecolor{amethyst}{rgb}{0.6, 0.4, 0.8}
\definecolor{burgundy}{rgb}{0.5, 0.0, 0.13}

\usepackage{lineno}

\title{Shapley values and machine learning to characterize metamaterials for seismic applications}

\author{\textbf{D.~\"Oniz}, 
\textbf{Y.~L.~Mo}, and \textbf{K.~B.~Nakshatrala} \\
Department of Civil and Environmental Engineering, University of Houston 77204.\\
\textbf{Correspondence to:} knakshatrala@uh.edu}

\begin{document}

\begin{titlepage}
  \drop=0.1\textheight
  \centering
  \vspace*{\baselineskip}
  \rule{\textwidth}{1.6pt}\vspace*{-\baselineskip}\vspace*{2pt}
  \rule{\textwidth}{0.4pt}\\[\baselineskip]
       {\Large \textbf{\color{burgundy}
       Shapley values and machine learning to characterize \\[0.1\baselineskip] metamaterials for seismic applications}}\\[0.3\baselineskip]
       \rule{\textwidth}{0.4pt}\vspace*{-\baselineskip}\vspace{3.2pt}
       \rule{\textwidth}{1.6pt}\\[\baselineskip]
       \scshape
       An e-print of the paper is available on arXiv. \par 
       \vspace*{0.5\baselineskip}
       Authored by \\[0.5\baselineskip]
           
  {\Large D.~\"Oniz\par}
  {\itshape Graduate student, Civil and Environmental Engineering, University of Houston, Texas 77204.}\\[0.3\baselineskip]  
 {\Large Y. L.~Mo\par}
{\itshape Professor, Civil and Environmental Engineering, University of Houston, Texas 77204.}\\[0.3\baselineskip] 
  {\Large K.~B.~Nakshatrala\par}
  {\itshape Department of Civil \& Environmental Engineering \\
  University of Houston, Houston, Texas 77204. \\ 
  \textbf{phone:} +1-713-743-4418, \textbf{e-mail:} knakshatrala@uh.edu \\
  \textbf{website:} http://www.cive.uh.edu/faculty/nakshatrala}\\
 
  \vfill
\begin{figure*}[htp]
\centering
    \includegraphics[height=3.0in]{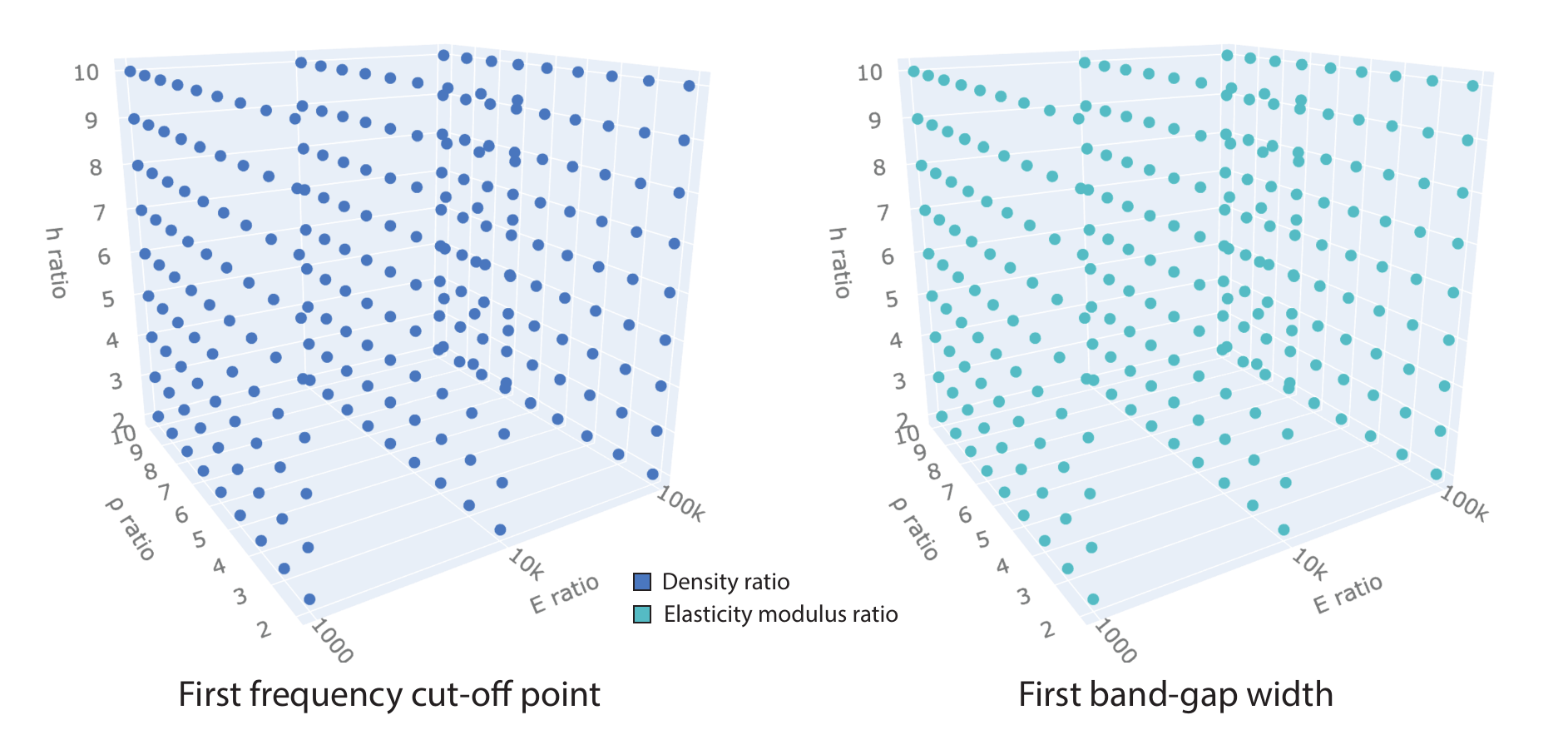}
\end{figure*}

\textit{
This figure illustrates the results from Shapley value analysis on a sonic metamaterial. For decreasing the first frequency cut-off point, increasing the density ratio is the only effective parameter. For increasing the width of the first band-gap,  the dominant parameter is the Young's modulus ratio.}
\vfill 

  {\small Computational \& Applied Mechanics Laboratory} \par
\end{titlepage}

\begin{abstract}
Given the damages from earthquakes, seismic isolation of critical infrastructure is vital to mitigate losses due to seismic events. A promising approach for seismic isolation systems is metamaterials-based wave barriers. Metamaterials---engineered composites---manipulate the propagation and attenuation of seismic waves. Borrowing ideas from phononic and sonic crystals, the central goal of a metamaterials-based wave barrier is to create band gaps that cover the frequencies of seismic waves. The two quantities of interest (QoIs) that characterize band-gaps are the first-frequency cutoff and the band-gap’s width. Researchers often use analytical (band-gap analysis), experimental (shake table tests), and statistical (global variance) approaches to tailor the QoIs. However, these approaches are expensive and compute-intensive. So, a pressing need exists for alternative easy-to-use methods to quantify the correlation between input (design) parameters and QoIs. To quantify such a correlation, in this paper, we will use Shapley values, a technique from the cooperative game theory. In addition, we will develop machine learning models that can predict the QoIs for a given set of input (material and geometrical) parameters.
\end{abstract}


\maketitle


\begin{tcolorbox}
\begin{center}
\textbf{SIGNIFICANCE STATEMENT}
\end{center}
\begin{itemize}
    \item The proposed sensitivity framework---based on Shapley values---requires minimal input data and is valuable for identifying dominant parameters in the design of seismic metamaterials, especially when the metamaterial response is not amenable to an analytical treatment (e.g., sonic crystals). 
    \item The paper also assesses various machine learning models in predicting quantities of interest (e.g., location of band-gaps) for seismic applications.
    \item The sensitivity framework and machine learning models will benefit a designer for tailoring a metamaterial to protect infrastructure against seismic disturbances.
\end{itemize}
\end{tcolorbox}

\setcounter{figure}{0}   
\section{INTRODUCTION AND MOTIVATION}
A crucial mission for urban planners and structural designers is to protect critical infrastructure (e.g., hospitals, bridges, power grids, to name a few) against earthquakes. It is also necessary to isolate sensitive equipment, such as MRI machines, in hospitals and laboratories from the ambient vibrations \citep{himmel2005structural}. The key strategy often employed to protect infrastructure and delicate instruments is through effective isolation systems that attenuate incoming disturbances, such as seismic waves and human-generated vibrations (e.g., from traffic). 

An emerging approach for seismic isolation is the use of metamaterials as wave barriers \citep{miniaci2016large,xiang2012periodic,huang2021experimental}. A metamaterial is an engineered composite that derives its functionalities through geometric layout of base materials \citep{cui2010metamaterials}. The advancement of photonic crystals for the light reflection led to the development of seismic metamaterials that can achieve band gaps, which will limit the propagation of the elastic waves within these gaps \citep{A.Khelif_2006}. Initially, phononic crystals were introduced in 1995 for filtering the Rayleigh waves \citep{brule2020emergence}. The structure of this material provides a full reflection of the elastic waves because of the periodically designed lattice \citep{R.H._Olsson_2011}. This periodic arrangement makes this metamaterial able to reflect the waves depending on whether its frequency creates constructive or destructive interference with the incoming waves. If the interference is destructive, the band-gap occurs. This phenomenon is commonly referred to as the Bragg scattering \citep{G.Yi_2016}. In 2000, \citet{liu2000locally} introduced the local resonance mechanism, which is also known as sonic crystals. For these crystals, the waves can be prevented by a structure consisting of spheres coated with a hard matrix. This mechanism differs from Bragg scattering as it does not require periodicity \citep{ungureanuphononic}. 

The basic idea of metamaterials-based wave barriers is to create band-gaps and ensure that these band gaps cover the frequencies of the incoming waves. It is essential to note two aspects: (a) Seismic waves have impulses in the low (5 – 50 Hz) and ultra-low (less than 5 Hz) frequency ranges. (b) The size of a unit cell is (inversely) proportional to the (frequency) wavelength. Thus, it will be practical to cover the seismic frequencies using lower band-gaps, preferably under the first band-gap, whereby maximizing the size of the wave barrier. This discussion implies that the first frequency cut-off (the lower boundary frequency between attenuation and propagation of the waves) and the width of the first band-gap (the distance between acoustic and optical branches of the waves) are the quantities of interest (QoIs).\footnote{Note that the first frequency cut-off is not the same as the first natural frequency of the structure.} The scientific question of importance about metamaterials-based wave barriers is: how do the input parameters (material and geometric properties) affect these quantities of interest? Said differently, which input parameters should we must vary to (i) lower the first frequency cut-off, and
(ii) increase the width of the first band-gap.

Many prior studies related to metamaterials-based seismic isolation systems have used analytical, experimental, statistical, or, more recently, machine learning approaches. Under an analytical study, the focus is to obtain analytical expressions for the dispersion relation---a mathematical relationship between frequency ($\omega$) and wave number ($k$) \citep{kittel1996introduction}. Solving the resulting dispersion relation provides the band diagram---a graph between $k$ and $\omega$ in the first Brillouin zone. The QoIs can then be determined from the band diagram. But except for canonical problems and simple layouts (e.g., layers), obtaining analytical dispersion relations is not possible; for example, metamaterials based on local resonance (i.e., sonic crystals) do not lend to analytical expressions.

Experimental studies on metamaterials have used either ultrasonic methods \citep{cheeke2017fundamentals} or shake table testing \citep{xiang2012periodic}. But both these approaches suffer from drawbacks. The former methods are unsuitable for seismic related studies, because of the limiting operating frequency ranges of the current ultrasonic instruments are covering the frequency ranges for ultrasonic waves only. There is no commercially available instrument capable of producing waves in the low and ultra-low frequency ranges. Next, a shake table test is expensive and labor-intensive. In addition, such a study needs special equipment (i.e., a shake table) and sophisticated instrumentation, available only in few universities and laboratories across the world (e.g., Pacific Earthquake Engineering Research, PEER, center \citep{ranf2001pacific}).

Under statistical methods, researchers have used parameter studies and global sensitivity analysis; both these methods have drawbacks. The na\"ive parametric approach---varying one input variable while keeping others constant---can not capture the interactions among the input variables. In contrast, global variance-based methods can account for such interactions \citep{saltelli2008global}. Recently, \citet{Witarto_2019} have used Sobol analysis---a popular method from the global analysis of variance---for phononic metamaterials to rank the input parameters affecting the QoIs. Although this approach accurately ranked the dominant parameters, a significant drawback was the need to perform numerous realizations ($\approx$100000) for the method to produce accurate results. So, this method is computationally prohibitive in carrying out those many realizations, especially for metamaterials based on sonic crystals. 

Lately, machine learning techniques have been employed to design metamaterials. For instance, \citet{wang2020deep} proposed a novel data-driven approach using deep generative modeling (a variational autoencoder) to achieve tailored microstructures for a metamaterial with the objective of achieving a target displacement profile. In another study, \citet{chen2022see} provided two machine-learning approaches to discover 2D metamaterials with user-defined frequency band gaps, offering logical rule-based conditions on unit cells; however, their designs are not suitable for seismic applications due to impractical fabrication feasibility and frequency ranges outside seismic requirements. Other notable machine-learning-related studies include \citep{kennedy2022machine,muhammad2022photonic,wu2020machine}. Nevertheless, these works focused on proposing general approaches to metamaterial design, considered frequency ranging that are not applicable for seismic applications, obtained material designs impractical to either fabricate in bulk and meet the design codes, or employed materials that are too expensive (e.g., gold). 

Thus, these drawbacks give rise to two specific needs: 
\begin{enumerate}[(a)]
    \item a simple procedure to quantify the effect of input parameters on the QoIs---the first frequency cut-off and the width of the first band-gap, and
    \item a robust method to predict QoIs in an economical and time-efficient way.
\end{enumerate} 
The aim of this paper is to address the said needs. Our \emph{innovation} is to bring ideas from cooperative game theory and machine learning. 

Our \emph{approach} to address the first need is to pose the relationship between the input variables and QoIs as a cooperative (mathematical) game. We then use the Shapley values---a technique from game theory---to calculate each input parameter’s contribution towards each QoI, thereby ranking the importance of the input parameters. Unlike the Sobol analysis, this approach does not need a large number of realizations when we have a small number of parameters \citep{narayanam2010shapley,jia2019towards}. To address the second need, we utilize machine learning algorithms to create a regression model for the QoIs. Once trained with sufficient data, the model can predict the QoIs for a different input data-set with a fraction of the time compared to the original analysis model. Bereft of these methods and knowledge, designing effective and economic metamaterials-based seismic isolation systems will remain unattainable.

The structure of the rest of this article is as follows. We will first provide preliminaries on phononic and sonic metamaterials and band-gap analysis of these metamaterials (\S\ref{Sec:S2_Shapley_preliminaries}). This discussion will be followed by the presentation of a sensitivity analysis framework using Shapley values (\S\ref{Sec:S3_Shapley_values}). The performance of the proposed framework will then be illustrated using numerical examples and exploratory data analysis for phononic crystals (\S\ref{Sec:S4_Shapley_NR_BS}) and sonic crystals (\S\ref{Sec:S5_Shapley_NR_LR}). After this, machine learning models will be presented to predict the QoIs (\S\ref{Sec:S5_Shapley_ML}). Finally, we will draw conclusions along with a discussion on possible future research directions (\S\ref{Sec:S6_Shapley_CR}).

\section{PRELIMINARIES: METAMATERIALS AND BAND-GAP ANALYSIS}
\label{Sec:S2_Shapley_preliminaries}

A metamaterial is a synthetic composite material system that is built to either imitate a response that is present in the nature or realize a new functionality \citep{engheta2006metamaterials}. Metamaterials are often designed to either achieve a distinct material property or tailor the response of the material system to an external stimuli. To provide a few examples of the former case, researchers have developed optical materials with negative refractive index \citep{smith2004metamaterials}, and materials with negative Poisson's ratio \citep{lakes1987foam,babaee20133d}. However, this paper addresses the other aspect of metamaterials---tailoring the response of a material system. 

The central theme of this paper is  manipulating the propagation of mechanical waves (such as seismic disturbances). Through an appropriate geometrical placement of the base materials, metamaterials can manipulate the propagation of waves through them. Depending on the material architecture, the waves can change direction, their intensity can diminish, or waves of certain frequencies will not be able to transmit at all: the resulting material system possesses band-gaps. For our purposes, a band-gap for a material system is a contiguous range of frequencies such that a wave within this frequency range will not be able to propagate through the system. In the parlance of condensed matter physics, it is the gap between acoustic and optical branches in the band diagram \citep{PierreA_Deymier_2012}.

\begin{figure}[h]
    \centering
        \includegraphics[height=1.3in]{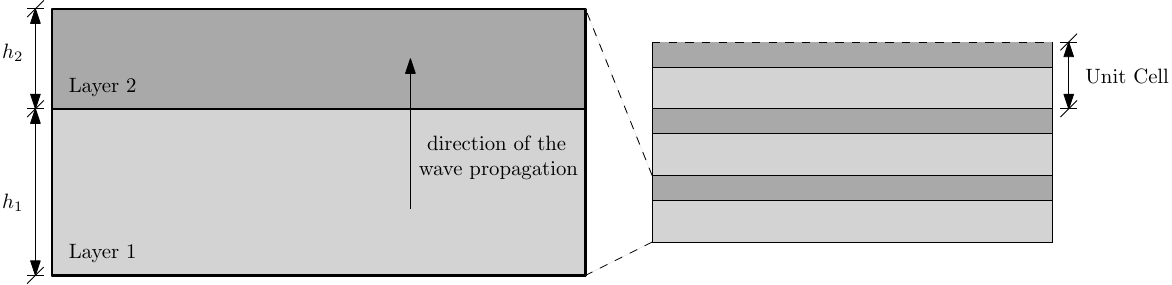}
        \caption{This figure shows the layered arrangement of base materials in a phononic metamaterial. The left figure shows a unit cell, which, in the case of a phononic crystal, is repeated infinite number of times. 
        \label{fig:layoutofphononic}}
\end{figure}

There are two main mechanisms for creating band-gaps: back (or Bragg) scattering and local resonance. Metamaterials that utilize back scattering are referred to as \emph{phononic} metamaterials, borrowing the name from phononic crystals \citep{S_Huber_2018}. Phononic crystals utilize periodic arrangement, as shown in \textbf{Fig.~\ref{fig:layoutofphononic}}. An incident wave upon such a material will undergo interference---a complex interaction among incident, transmitted and reflected components of the wave. This interference can be either destructive or constructive \citep{S_Huber_2018}. If the interference is destructive, the material can neutralize the amplitude of the wave, thereby creating  frequency band-gaps \citep{PierreA_Deymier_2012}. An attractive feature of a phononic metamaterial is its simple layout, which allows an analytical treatment of dispersion analysis. Using the transfer matrix method, \citet{Witarto_2019} obtained the following dispersion equation---a mathematical relation or graph between the wave number $k$ and (angular) frequency $\omega$---for a layered phononic metamaterial: 
\begin{equation}
    \cos(k h) = \cos\left(\frac{\omega h_1}{C_1}\right)\cos\left(\frac{\omega  h_2}{C_2}\right) - \frac{1}{2}\left(\frac{\rho_1 C_1}{\rho_2 C_2}+\frac{\rho_2 C_2}{\rho_1 C_1}\right)\sin\left(\frac{\omega h_1}{C_1}\right)\sin\left(\frac{\omega h_2}{C_2}\right)
    \label{equation:dispersionbraggs}
\end{equation}
where $h$ represent the thickness of a unit cell, while $h_n$, $\rho_{n}$, and $C_{n}$, respectively, denote the thickness, density, and wave speed in the $n$th layer. 

\begin{figure}[h!]
    \includegraphics[width=0.7\textwidth]{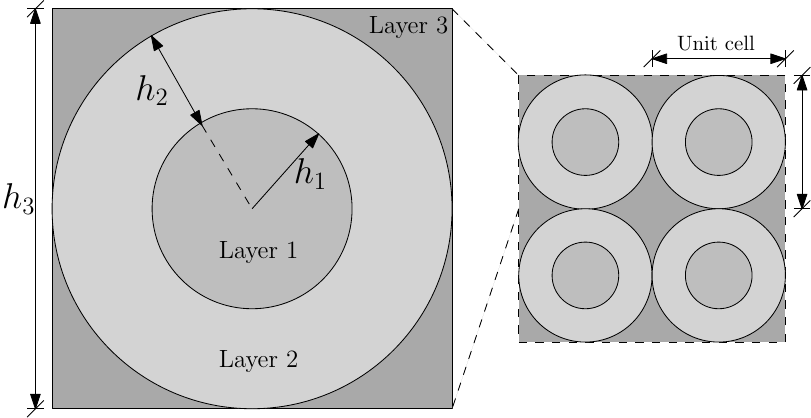}
    \caption{This figure shows the arrangement of base materials in a sonic crystal. The left depicts a unit cell, which is repeated to realize sonic crystal in two dimensions. \label{Fig:Shapley_Sonic_crystal}}
\end{figure}

On the other hand, \emph{sonic} metamaterials utilize local resonance to achieve band-gaps. In a seminal paper entitled: “Locally Resonant Sonic Materials,” appeared in the journal of Science, the concept of local resonance is introduced as a way to realize band-gaps. Their key idea hinges on constructing a (sonic) crystal with a lattice constant equal to two order of magnitudes smaller than the wavelength \citep{liu2000locally}. The resulting sonic metamaterial comprised three parts. The core is made of lead which is covered with a thin layer of silicone rubber; this layer is commonly referred to as the resonator. The core and thin layer are put into a cube of epoxy, which has a stiffness between the lead and silicone rubber \citep{liu2000locally}. The layout of such a sonic crystal is shown in  \textbf{Fig.~\ref{Fig:Shapley_Sonic_crystal}}. When the resonance frequency of the resonator and the propagating seismic waves interact, a band-gap is created. The locations and widths off the frequency band-gap occurs are manipulated by tuning the parameters of the resonating (thin layer), elastic (epoxy) and inertial (core) components of the sonic metamaterial \citep{raghavan2013local}.

As mentioned above, we already know an analytical expression for dispersion relation for layered media. We therefore use this information to establish the utility and accuracy of Shapley value. After that, the proposed technique will be applied to the sonic metamaterial -- for which an analytical dispersion relation does not exist.

\section{A SENSITIVITY ANALYSIS FRAMEWORK BASED ON SHAPLEY VALUES}
\label{Sec:S3_Shapley_values}
Game theory is a mathematical approach often utilized in decision-making, especially in situations involving social interactions. This approach is widely used in economics, politics, the strategic game industry, and fields such as engineering and computer sciences \citep{Bhuiyan_2018}. In the language of game theory, a situation is referred to as a `game' and each party involved in the situation is referred to as a `player'  \citep{R_Giles_2010}. The study of game theory is broadly classified into cooperative and non-cooperative game theory. In a cooperative game, each player agrees to work together to achieve the desired result. On the other hand, non-cooperative game theory considers situations in which each player wants to achieve success, regardless of the other players decisions; that is, there can be competition among players. For further details on game theory, see \citep{osborne1994course}. 

In a seminal paper, \citet{shapley1953value} proposed an approach to evaluate a game by calculating the marginal contribution of each player; these marginal contributions are quantified by Shapley values \citep{S_Hart_1989}. Shapley value is helpful to determine the value of each player when there is cooperation among them \citep{J_Castro_D_Gamez_2009}. When this idea is integrated into the field of metamaterials, a game---where the material properties are the players---can be constructed. 

In this paper, cooperative game theory is utilized to characterize the mechanics of metamaterials, and the Shapley values will provide the relative importance of the input parameters over the design. Each of the two QoIs—decreasing the first-frequency cut-off, and increasing the first band-gap—is posed as a separate game. The input parameters are the players of the game. For each game, players’ pay-off will be calculated which provide the rankings of the parameters according to their Shapley values.

Since game theory is foreign to the field of mechanics, we will provide several examples to illustrate the main concepts and modifications used in this paper.

\subsection{Illustrative examples}
The first example provides an overview of cooperative game theory: introduces the notation such as a game and players. The second example demonstrates the Shapley value analysis, while the third shows the importance of satisfying the super-additive property. In mathematics, a function $v(\cdot)$ is said to satisfy the super-additive property if 
\begin{align}
    v(A \cup B) \geq v(A) + v(B) 
\end{align}
where $A$ and $B$ are sets. In game theory, $v(\cdot)$ denotes the characteristic function of the game. All the three games have a similar scenario but varies by minor explanatory differences. 

\subsubsection{Illustrating a cooperative game.} Three students---Alice (A), Bill (B) and Charlie (C)---have to write together a report of at least 70 pages within a week for a competition. Individually, these students can write 20, 27 and 35 pages within a week, respectively. If they work in groups of two:  Alice and Bill can write 55 pages, Alice and Charlie can write 62 pages, and Bill and Charlie together can write 74 pages. All three can write together and accomplish 100 pages in a week. Then the decision to make is: what is the most economical way to meet the requirements for the competition---a report of at least 70 pages within a week. 

\textbf{Table} \ref{Table:Shapley_writing_a_report_game} poses the mentioned scenario as a mathematical \emph{game} of three \emph{players}. Two combinations can meet the competition's requirements---as given in bold: they all can work together to produce a 100-page report, or only Bill and Charlie can work together to produce a 74-page report. Thus, the economical option---to minimize labor---is for Bill and Charlie to work in a group of two.

This example perfectly illustrates the basic concept of cooperative game theory. In this example, number of written pages has no utility, where the only important thing is to qualify for the contest by writing at least 70 pages. 

\begin{table}[h]
\caption{A cooperative game of \emph{writing a report}.\label{Table:Shapley_writing_a_report_game}}
\begin{tabular}{c|c|c|c|c|c|c|c}
\textbf{Combination} & A & B & C & A\&B & A\&C & B\&C & A\&B\&C\\
\hline
\textbf{Pages written} & 20 & 27 & 35 & 55 & 62 & \textbf{74} & \textbf{100}
\end{tabular}
\end{table}

\subsubsection{Illustrating Shapley value analysis.} An additional requirement for the competition was defined later as minimum three participants are required per group. Therefore, even though it is not the most economical combination, Alice, Bill and Charlie joined the competition all together. In the end, they won the competition and got the prize of 1000 dollars. They want to split the prize depending on the individual contribution of each. This contribution will be determined only by the number of pages written by that individual, not by the importance of what was written. The Shapley value analysis is shown in \textbf{Table} \ref{Table:Shapley_prize_sharing_Shapley_values}. In this table, the writing order A-B-C represents a scenario where Alice writes her share first (20 pages), followed by Bill until their combined page number (55 pages) is reached. Charlie writes last until 100-page report is done. After calculating the Shapley value---average marginal contribution---of each, multiplying the dominance percentage with the prize---1000 dollars---will give the prize each should take.

\begin{table}[h]
\caption{Shapley value analysis for the cooperative game of \emph{sharing the \$1000 prize}. 
\label{Table:Shapley_prize_sharing_Shapley_values}}
\centering
\begin{tabular}{c|c|c|c}
\textbf{Writing order} &\textbf{A} & \textbf{B} & \textbf{C}\\
\hline
A-B-C & 20 & 35 & 45  \\
A-C-B & 20 & 38 & 42 \\
B-A-C & 28 & 27 & 45 \\
B-C-A & 26 & 27 & 47 \\
C-A-B & 27 & 38 & 35 \\
C-B-A & 26 & 39 & 35 \\ \hline
\textbf{Total} & 147 & 204 & 249\\ \hline
\textbf{Shapley value (in pages)} & 24.5 & 34 & 41.5\\ \hline
\textbf{Dominance (in \%)} & 24.5 & 34 & 41.5\\ \hline
\textbf{Money award (in \$)} & 245 & 340 & 415\\
\end{tabular}
\end{table}

\subsubsection{Illustrating super-additive property}
This example considers human interactions. Just as in the original scenario, the students individually can write 20, 27, and 35 pages in a week, respectively. However, Alice does not work well in a group---she slows down the process by effecting others. When she pairs with Bill, both can write only write 10 pages in a week. Alice and Charlie together can write 17 in a week. On the other hand, Bill and Charlie work great together and can write up to 74 pages within a week. When all three students work together, they can complete 70 pages in a week. \textbf{Table} \ref{Table:Shapley_writing_a_report_game_super_additive} summarizes this scenario. 

In this scenario, when Alice is involved into a combination, it does not satisfy super-additive property.

\begin{table}[h]
\caption{A cooperative game of \emph{writing a report without satisfying super-additive property}.\label{Table:Shapley_writing_a_report_game_super_additive}}
\begin{tabular}{c|c|c|c|c|c|c|c}
\textbf{Combination} & A & B & C & A\&B & A\&C & B\&C & A\&B\&C\\
\hline
\textbf{Pages written} & 20 & 27 & 35 & 10 & 17 & \textbf{74} & \textbf{70}
\end{tabular}
\end{table}

For this competition, they have two ways to achieve their goal---working all together or working without Alice. If they decided to work all together and won the same prize, they again decide to split it according to their Shapley values. The Shapley analysis is given in \textbf{Table} \ref{Table:Shapley_prize_sharing_Shapley_values_super_additive}.

\begin{table}[h]
\caption{Shapley value analysis for the cooperative game of \emph{sharing the \$1000 prize without satisfying super-additive property}. \label{Table:Shapley_prize_sharing_Shapley_values_super_additive}}
\begin{tabular}{c|c|c|c}
\textbf{Writing order} &\textbf{A} & \textbf{B} & \textbf{C}\\
\hline
A-B-C & 20 & -10 & 60  \\
A-C-B & 20 & 53 & -3 \\
B-A-C & -17 & 27 & 60 \\
B-C-A & -4 & 27 & 47 \\
C-A-B & -18 & 53 & 35 \\
C-B-A & -4 & 39 & 35 \\ \hline
\textbf{Total} & -3 & 189 & 234\\ \hline
\textbf{Shapley value (in pages)} & -0.5 & 31.5 & 39\\ \hline
\textbf{Dominance (in \%)} & -0.71 & 45 & 55.71\\ \hline
\textbf{Money award (in \$)} & -7.1 & 450 & 557.1\\
\end{tabular}
\end{table}

As it can be seen from the Table \ref{Table:Shapley_prize_sharing_Shapley_values_super_additive}, even if Alice contributed and they won, she has to give back 7.1 dollars to others and get no prize after all. In addition to this, the contribution of Alice is negative. Therefore, it makes it difficult to rank the importance of these students as they are effecting the result the opposite way. There is a chance that Alice may be effecting it in a more negative way than Charlie does in a positive way. 

When super-additive property is not satisfied, data can be modified to get clearer results. Modifying the data-set can help students to determine who contributed the most and how they should split the prize. When combination of A and B is less than either A or B, the maximum value of these two should be taken as the new combination. \textbf{Table} \ref{Table:Shapley_writing_a_report_game_modified} explains the modifying method, while \textbf{Table} \ref{Table:Shapley_prize_sharing_Shapley_values_modified} shows the Shapley value analysis with the modified data. 

\begin{table}[h]
\caption{A cooperative game of \emph{writing a report with modified data-set}.\label{Table:Shapley_writing_a_report_game_modified}}
\begin{tabular}{c|c|c|c|c|c|c|c}
\textbf{Combination} & A & B & C & A\&B & A\&C & B\&C & A\&B\&C\\
\hline
\textbf{Pages written} & 20 & 27 & 35 & 27 & 35 & \textbf{74} & \textbf{74}
\end{tabular}
\end{table}

\begin{table}[h]
\caption{Shapley value analysis for the cooperative game of \emph{sharing the \$1000 prize with modified data-set}. \label{Table:Shapley_prize_sharing_Shapley_values_modified}}
\begin{tabular}{c|c|c|c}
\textbf{Writing order} &\textbf{A} & \textbf{B} & \textbf{C}\\
\hline
A-B-C & 20 & 7 & 47  \\
A-C-B & 20 & 39 & 15 \\
B-A-C & 0 & 27 & 47 \\
B-C-A & 0 & 27 & 47 \\
C-A-B & 0 & 39 & 35 \\
C-B-A & 0 & 39 & 35 \\ \hline
\textbf{Total} & 40 & 178 & 226 \\ \hline
\textbf{Shapley value (in pages)} & 6.67 & 29.67 & 37.67 \\ \hline
\textbf{Dominance (in \%)} & 9.01 & 40.09 & 50.91 \\ \hline
\textbf{Money award (in \$)} & 90 & 401 & 509 \\
\end{tabular}
\end{table}

Modifying the calculations---mainly focused on Alice's contribution to achieve the goal instead of her negative effect---resulted in her having a positive effect. Therefore, she is able to get some of the award for her work or at least making the team qualify for the contest. However, as the focus is on the positive effect for this game, when compared with others her contribution is considered as less valuable. 

For the case where the qualification is changed to at least two people for the competition, Bill and Charlie chose to work together without Alice. By doing so they can share the award based on their contribution. \textbf{Table} \ref{Table:Shapley_prize_sharing_Shapley_values_BandC} below explains the distribution of the award between Bill and Charlie. As it can be seen from the table, by neglecting Alice, Bill and Charlie can write more pages and gain much more money from the award. Therefore, if applicable, in some cases working with less parameters can lead to more economical results. 

\begin{table}[htp]
\caption{Shapley value analysis for the cooperative game of \emph{sharing the \$1000 prize without Alice}. \label{Table:Shapley_prize_sharing_Shapley_values_BandC}}
\begin{tabular}{c|c|c}
\textbf{Writing order} & \textbf{B} & \textbf{C}\\
\hline
B-C & 27 & 47 \\
C-B & 39 & 35 \\ \hline
\textbf{Total} & 66 & 82\\ \hline
\textbf{Shapley value (in pages)} & 33 & 41\\ \hline
\textbf{Dominance (in \%)} & 44.59 & 55.41\\ \hline
\textbf{Money award (in \$)} & 446 & 554\\
\end{tabular}
\end{table}

In conclusion, when the super-additive property is not satisfied, there are two possible solutions to overcome this problem: (i) modifying the data-set, or (ii) neglecting the negative parameters. When a parameter has no positive effect both individually or within a combination, modification method used in this paper fully neglects that parameter.



\subsection{Application on a simple continuous model}
Till now, we have applied the Shapley value analysis to discrete problems. We now consider a continuous model and demonstrate how to conduct Shapley value analysis on such a problem. 

The mathematical model is given by
\begin{align}
\label{Eqn:Shapley_continuous}
    f(x_1,x_2) = 3 x_1^2 + x_2^2 - x_1 x_2, \quad 
    x_1, x_2 \in [0,10]
\end{align} 
The range for the two input parameters $x_1$ and $x_2$ is set to $[0, 10]$. We lay out a grid with a spacing of $0.1 \times 0.1$. At each grid point, we perform Shapley value analysis to quantify the relative importance of the input parameters, expressed as percentages.

\begin{figure}[htp]
    \includegraphics[width=0.6\textwidth]{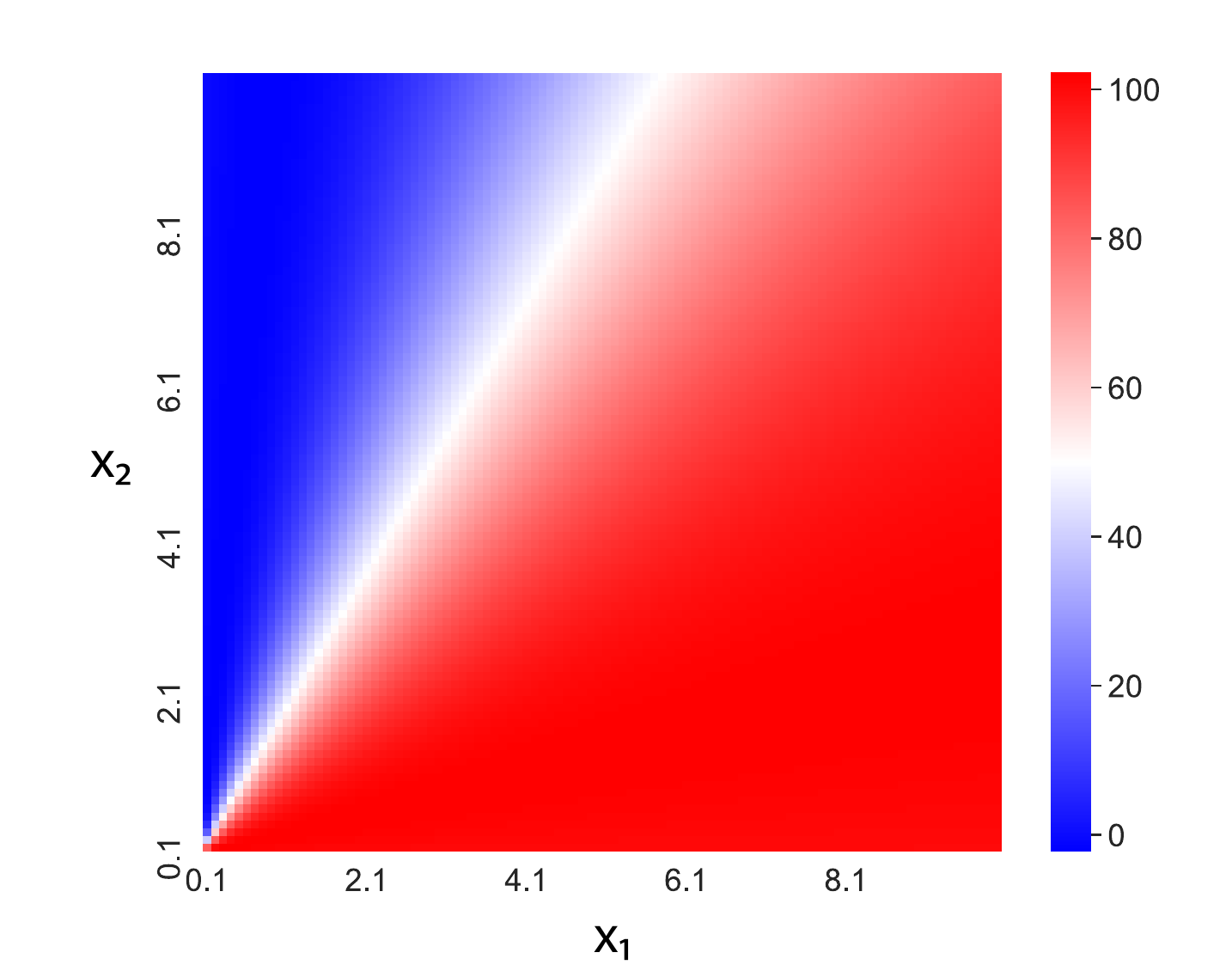}
    \caption{This figure presents the Shapley value analysis results for the mathematical model \eqref{Eqn:Shapley_continuous}. The graph plots the $x_1$ dominance in percentage. Thus, in the regions with lower percentages [$\leq 50$], the parameter $x_2$ is dominant.\label{fig:math_model}}
\end{figure}

\textbf{Figure \ref{fig:math_model}} illustrates the dominance of $x_1$ by plotting its Shapley value as a percentage. Regions where $x_1$ contributes more than 50\% are depicted in red, indicating its dominant influence. Conversely, blue regions represent where $x_2$ is dominant, with white indicating the transition between these influences. The predictions from the Shapley value analysis align with the analytical solution. Moreover, Shapley value analysis can be applied even in scenarios lacking analytical solutions, such as designing seismic wave barriers using metamaterials. In such cases, graphs akin to Fig.~\ref{fig:math_model} offer valuable insights to design engineers navigating the design space.

For instance, \textbf{Eq.~\eqref{Eqn:Shapley_continuous}} represents a material design involving two materials, where $x_1$ denotes the density ratio and $x_2$ represents the thickness ratio between these materials. In practical scenarios, such as when specific types of steel and rubber are available, adjusting the density ratio may necessitate purchasing new materials, which is economically impractical. By referencing the information depicted in Fig.~\ref{fig:math_model}, we can identify the ranges where modifying the thickness ratio effectively influences the design. Suppose the available density ratio is 4. Maintaining the thickness ratio above 7.5 grants us the advantage of having thickness as the dominant parameter, allowing us to modify the design as needed. In conclusion, identifying the dominant parameter through the Shapley value helps us determine how to adjust the design characteristics as desired, in an economical and time-efficient manner.

\section{EXPLORATORY DATA ANALYSIS (EDA) FOR BRAGG  SCATTERING}
\label{Sec:S4_Shapley_NR_BS} 

\subsection{Material properties and design of the phononic crystal}

The material under investigation is a one-dimensional metamaterial comprised of two base materials. Figure~\ref{fig
} illustrates the arrangement of the phononic crystal. The effective design parameters considered include the modulus of elasticity, density, and thickness ratios. These ratios are denoted as $x_{\text{layer1}}:x_{\text{layer2}}$, with $x$ representing the characteristic of each layer—Young's modulus, density, and thickness. Poisson's ratio remains constant, as Sobol analysis conducted by \citet{Witarto_2019} has demonstrated its negligible impact on the design.

\begin{figure}[htp]
    \hspace{0.5in}\includegraphics[width=0.75\textwidth]{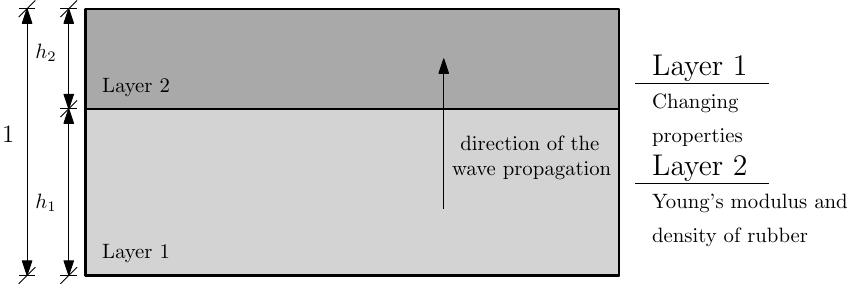}
    \caption{This figure shows the arrangement of the base materials in the phononic metamaterial for the exploratory data analysis. The base material for layer 2 is rubber; only the geometrical properties of layer 2 are changed, while (dimensionless) unit cell width is fixed at 1. The properties of layer 1 are varied by altering the ratios. The changing parameters for the design are Young's modulus ratio, density ratio, and thickness ratio. These ratios are calculated as $x_{layer 1}:x_{layer 2}$ for a given property $x$.
    \label{fig:arrangement_braggs}}
\end{figure}

In this metamaterial, the properties of the second material, embodied by layer 2, are configured to mirror those of rubber, with its Young's modulus set at 3.49 MPa and its density at 1100 $\mathrm{kg/m^3}$. Throughout the analysis, only the geometrical properties of the second material undergo alteration. However, the dimensionless thickness of the unit cell, denoted as $h_1 + h_2$, remains fixed at 1. Conversely, for the first material, represented by layer 1, properties are systematically varied to explore their impact on the design. The baseline value for all parameter ratios in the Shapley value analysis is established at 0.1.

For the design, the primary objectives include achieving a lower first-frequency cut-off point to encompass ultra-low frequency ranges and a broader band-gap width to accommodate a wider range of seismic wave frequencies. Hence, when specific parameters are labeled as negative in this paper, it indicates that these parameters lead to higher first-frequency cut-off points or narrower band-gap widths, which are undesirable outcomes.

\subsection{Data set}
The data used in this study was generated using a Python code. This code first applied dispersion relation calculations to determine the quantities of interest (QoIs) and then applied Shapley value analysis to identify the dominant parameter among the three input ratios. The two QoIs are decreasing the first frequency cut-off point and increasing the first band-gap. For this analysis, the dataset consisted of ranges of the parameters, as detailed in \textbf{Table \ref{Table:dataset_bragg_range}}.

\begin{table}[h]
\caption{Parameter ranges for Shapley value analysis of \emph{phononic crystals}. \label{Table:dataset_bragg_range}}
\begin{tabular}{l|r}
\textbf{Parameter} & \textbf{Range}\\
\hline
Young's modulus ratio & 0.1--50000\\
\hline
Density ratio &  0.1--9.5\\
\hline
Thickness ratio & 0.1--11\\
\end{tabular}
\end{table}

\subsubsection{Visualizing labeled data for first frequency cut-off point.}

\begin{figure}[ht!]
    \includegraphics[width=1.0\textwidth]{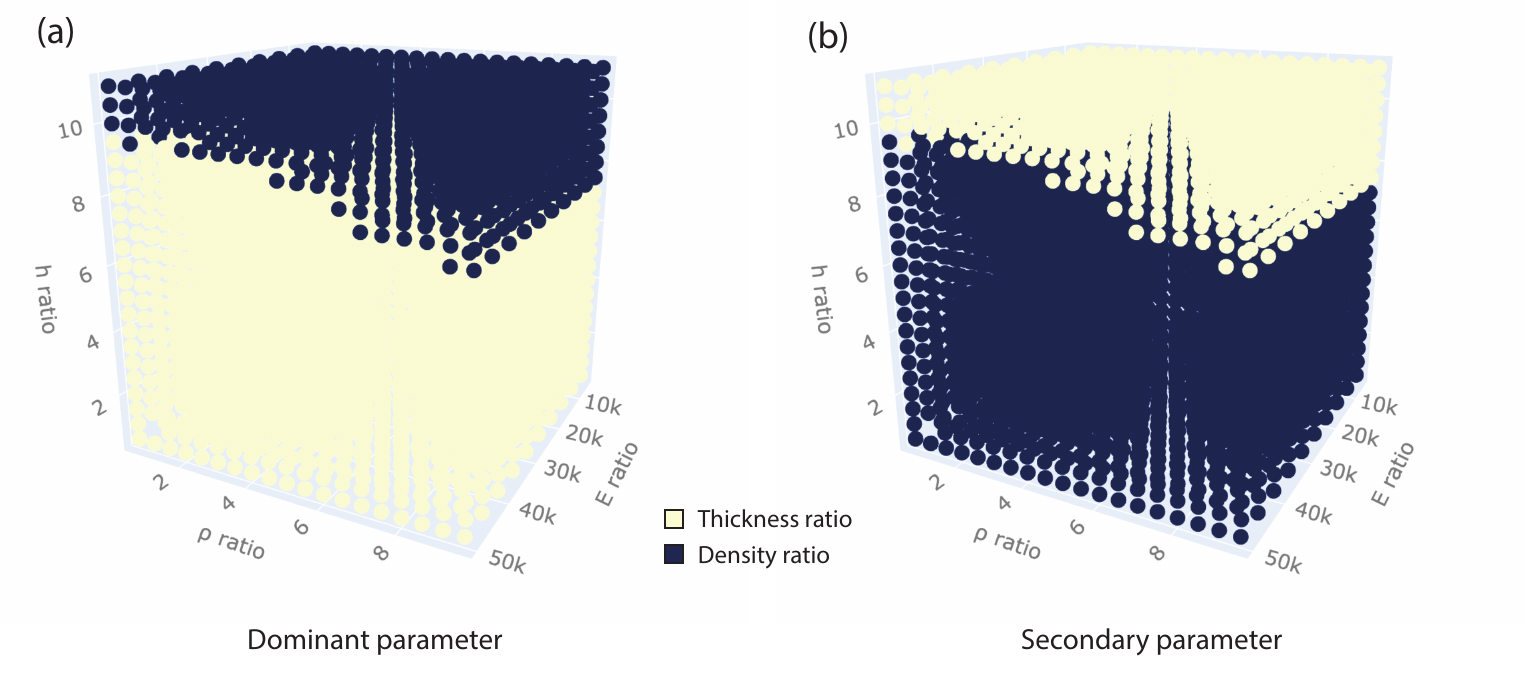}
   \caption{\textsf{Shapley value analysis for a phononic metamaterial.} This figure presents the parameter ranking for the \emph{first} QoI---decreasing the first frequency cut-off point. The dominant and secondary parameters are shown in subfigures (a) and (b), respectively. For this QoI, elasticity modulus has no effect over the design. For the lower ranges of thickness ratio [$\leq 9.5$] and density ratio [$\leq 2$], thickness ratio governs the design. Above this range, the density ratio is the dominant parameter. However, this dominance zone of thickness ratio reduces as the density ratio increases.}
   \label{Fig:dominantffcop}
\end{figure}

\textbf{Figure~\ref{Fig:dominantffcop}} presents, utilizing 3D colormaps, the results obtained from the Shapley value analysis for the first QoI, which aims to decrease the first frequency cut-off point. Parameters are characterized in terms of ratios.

Since the data modification omitted a parameter that showed no positive effect, the marginal contribution for Young's modulus was obtained as zero within the specified ranges. Consequently, there was no tertiary parameter influencing the design. Said differently, increasing Young's modulus, either alone or in combination with other parameters, only raises the first frequency cut-off point rather than lowering it. The dominant and secondary parameters of the design are illustrated in Fig.~\ref{Fig:dominantffcop}. The absence of data when both density and thickness ratios are equal to 1 is attributed to the lack of band-gap generation. In these diagrams, yellow and blue denote the thickness and density ratios, respectively. In the lower ranges of thickness ratio [$\leq 9.5$] and density ratio [$\leq 2$], the thickness ratio emerges as the dominant parameter for reducing the first frequency cut-off point. Beyond these thresholds, the density ratio begins to dictate the design, extending its dominance range as it increases.

\subsubsection{Visualizing labeled data for the width of first band-gap.} 

\begin{figure}[ht!]
    \includegraphics[width=1.0\textwidth]{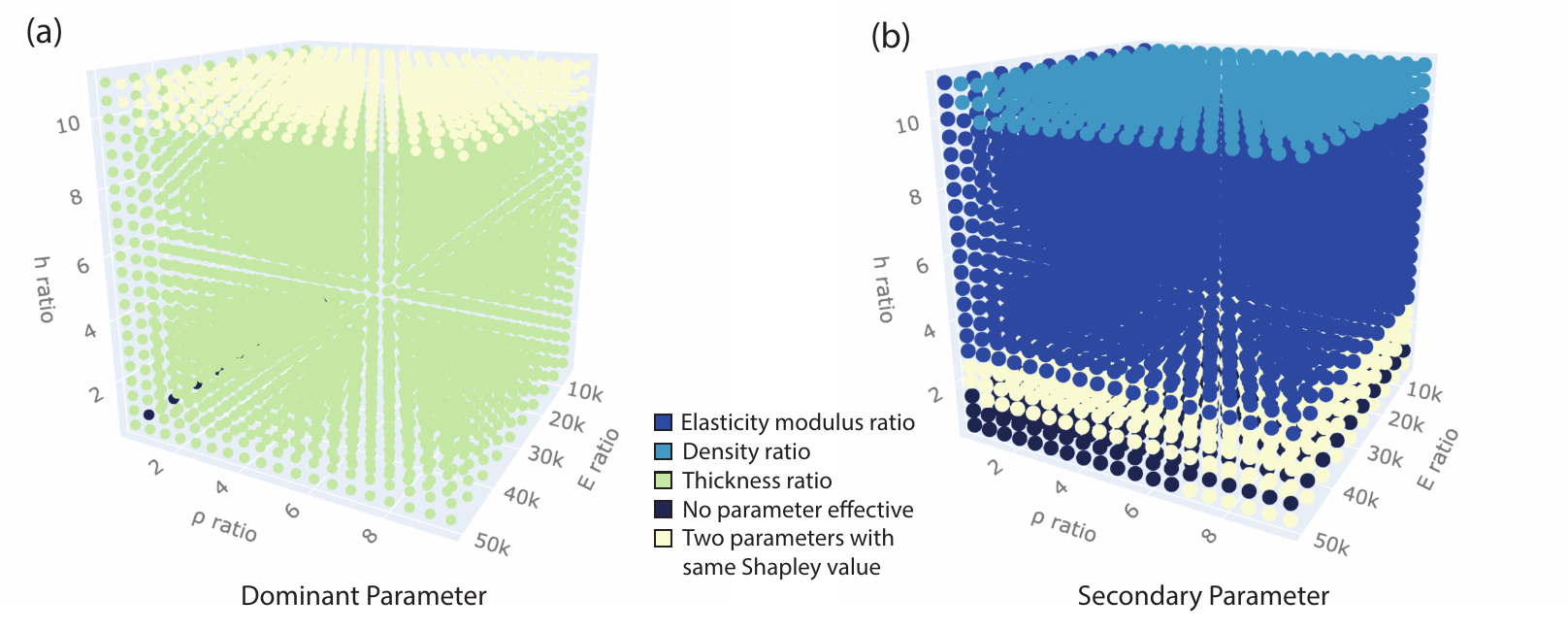}
    \caption{\textsf{Shapley value analysis for a phononic metamaterial.} This figure illustrates the ranking of parameters for the \emph{second} QoI---increasing the band-gap width. For phononic crystal, no band-gap data is obtained when both density and thickness ratios are equal to 1. For the dominant parameter (subfigure a), for the lower ranges of thickness ratio [$\leq 9.5$], thickness ratio governs the design. Above these ranges, thickness and Young's modulus ratios have the same dominance over the design. For the secondary parameter (subfigure b), when thickness ratio is equal to 1, the other two parameters become ineffective. However, for lower ranges of thickness ratio [$\leq 2.5$], density and elasticity modulus ratios have the same effect over the design. For intermediate ranges [2.5--10], elasticity modulus is the second dominant parameter. Above this point, density becomes the secondary parameter, as other two parameters are both dominant with same Shapley value.}
    \label{fig:dominantsecondfbg}
\end{figure}

\textbf{Figure~\ref{fig:dominantsecondfbg}} depicts, through a colormap, the dominant and secondary parameters influencing the increment of the first band-gap width. Parameters in all the graphs are represented in terms of ratios. 

In the colormap representing parameter rankings, each color corresponds to a distinct parameter. Dark blue indicates no dominant parameter within the range, while the modulus of elasticity ratio, density ratio, and thickness ratio are depicted in blue, light blue, and green, respectively. Yellow is utilized to denote ranges where two parameters share the same Shapley value. As illustrated in Fig.~\ref{fig:dominantsecondfbg}\textcolor{blue}{a}, the thickness ratio emerges as the dominant parameter governing the design within lower thickness ratio ranges [$\leq 9.5$]. Beyond this range, both the elasticity modulus and thickness ratios make an equal marginal contribution to the design. The absence of a dominant parameter when both density and thickness ratios are equal to 1 stems from the lack of band-gap generation. As depicted in Fig. \ref{fig:dominantsecondfbg}\textcolor{blue}{b}, no secondary parameter is evident, both when the thickness ratio is equal to or less than 1 and within lower ranges of density ratio [$\leq 6.5$]. Beyond these ranges, both density and elasticity modulus ratios exert equal influence on the design until the thickness ratio reaches 3. Beyond this point, Young's modulus begins to take control (blue color in Figs.~\ref{fig:dominantsecondfbg}\textcolor{blue}{a} and \ref{fig:dominantsecondfbg}\textcolor{blue}{b}).

\subsection{Discussion on EDA for Bragg Scattering.}

In this section, we utilized both the dispersion relation and Shapley value analysis to comprehend the behavior of the phononic metamaterial. Initially, we determined the band-gap characteristics of the crystal by applying the dispersion relation across various ranges and combinations of design parameters. Subsequently, Shapley value analysis was applied individually to each QoI to pinpoint the dominant parameter capable of affecting the desired design alterations.

The main findings from this exploratory data analysis of the \emph{Bragg scattering} dataset are outlined as follows: In terms of reducing the \emph{first frequency cut-off point}, the Young's modulus ratio exhibits no influence on the dataset. For density ranges at or below 8.5 and thickness ratios at or below 9.5, the predominant factor is the thickness ratio. However, as these parameters increase, the density ratio begins to dictate the design. Consequently, to attain a lower first frequency cut-off point, the following strategies are recommended: (i) Maintain the Young's modulus ratio at the same level or within lower ranges, as feasible within the design constraints. (ii) When working with lower density and thickness ratios, increasing the thickness ratio can result in a lower first frequency cut-off point. (iii) For higher ranges of both thickness and density ratios, increasing the density ratio yields superior outcomes.
~
To \emph{increase the width of first band-gap}, the thickness ratio emerges as the sole effective parameter within lower ranges of both density ($\leq 6.5$) and thickness ratios ($0.1-1$). Within intermediate ranges ($1.2-9.5$), the thickness ratio remains the dominant factor. However, as the thickness ratio surpasses 9.5, both Young's modulus and thickness ratios exert equal influence over the design. Therefore, to attain a larger band-gap, the following strategies are recommended: (i) When operating within the lower ranges of thickness ($0.1-1$) and density ratios ($0.1-6.5$), a broader band-gap can be achieved by solely increasing the thickness of layer 1. Any adjustments to the other two parameters within these ranges will only diminish the width of the band-gap. (ii) Increasing the thickness ratio will also expand the band-gap for higher ranges, up to the point where the thickness ratio exceeds 9.5. Beyond this threshold, the Young's modulus ratio can also be increased.

\subsubsection{Comparison of results with Sobol analysis}

In previous research, \citet{Witarto_2019} explored the primary and influential parameters affecting the design of phononic crystals using Sobol analysis within the specified ranges outlined in Table \ref{Table:sobol_range}. Notably, Poisson's ratio range from the prior study is omitted here, as the parameter did not affect QoIs; hence, it was not considered in this paper. The parameter ranges utilized for the Sobol analysis align with our observational range. However, the range employed for Shapley value analysis exhibits a higher upper limit for Young's modulus and a lower limit for density ratio. Moreover, a significant distinction between Sobol and Shapley value analysis lies in Sobol's ability to capture parameter interactions, a facet not accounted for in Shapley value analysis.

\begin{table}[h]
\centering
\caption{Parameter ranges for Sobol analysis from \citet{Witarto_2019}. \label{Table:sobol_range}}
\begin{tabular}{l|r}
\centering
\textbf{Parameter} & \textbf{Range}\\
\hline
Young's modulus ratio & 10--10000\\
\hline
Density ratio &  1--1000\\
\hline
Thickness ratio & 0.11--9\\
\end{tabular}
\end{table}

After conducting the Sobol analysis, \citet{Witarto_2019} reported the following findings: For the \emph{first frequency cut-off point}, the predominant parameter identified was the density ratio, followed by the interaction between density and thickness ratio, and then thickness ratio alone. Moreover, for P-waves, Poisson's ratio and its interaction with other parameters were found to significantly influence the design. Regarding the \emph{width of the first band-gap}, the dominant parameter was determined to be the thickness ratio. Following this, the interaction between Young's modulus and thickness ratio emerged as a secondary factor. The tertiary parameter was the interaction between density and thickness ratios, with subsequent influence from Young's modulus ratio and density ratio, individually. Furthermore, for P-waves, Poisson's ratio was also observed to have a notable effect on the design.

The comparison between the results obtained from Sobol analysis and our findings on Shapley value analysis can be outlined as follows. For the \emph{first frequency cut-off point}:
\begin{enumerate}
\item[(i)] Shapley value analysis identifies density ratio as the dominant parameter for higher ranges. Given that Sobol analysis extends beyond our density ratio range, it suggests that density likely remains influential across most values within the range. Thus, the prominence of density as a dominant parameter in Sobol analysis can be justified, as the sensitivity analysis does not discern transitions between dominant parameters.
\item[(ii)] Additionally, Shapley value analysis indicates that Young's modulus ratio does not impact the reduction of the first frequency cut-off point, a finding mirrored by Sobol analysis, where Young's modulus ratio does not emerge as an effective parameter.
\end{enumerate}
For the \emph{width of the first band-gap}:
\begin{enumerate}
\item[(i)] Sobol analysis identifies thickness ratio as the dominant factor. Considering the observation range of the Sobol analysis, it consistently aligns with the ranges where the thickness ratio governs the design in the Shapley value results.
\item[(ii)] Notably, when interactions are not factored in, secondary and tertiary parameters predicted by Sobol analysis align closely with the results of Shapley value analysis.
\end{enumerate} 

\section{EXPLORATORY DATA ANALYSIS (EDA) FOR LOCAL RESONANCE}
\label{Sec:S5_Shapley_NR_LR} 

\subsection{Material properties and design of the sonic crystal}
\begin{figure}[htp!]
    \includegraphics[width=0.6\textwidth]{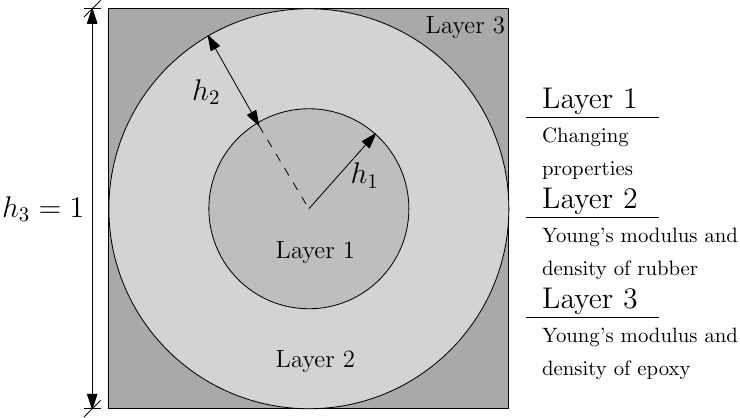}
    \caption{This figure shows the arrangement of base materials.   During the Shapley value analysis, only the properties of layer 1 are varied. The design parameters are taken as $x_{\text{layer1}}:x_{\text{layer2}}$, where $x$ is the property of the material---such as density, Young's modulus, and thickness.}
    \label{fig:lr}
\end{figure}

We considered a 2D sonic metamaterial comprising three base materials. \textbf{Figure~\ref{fig:lr}} shows the arrangement of base materials.

The base material in the second layer is taken to be rubber. Since the material is fixed in this layer, only the geometrical properties (e.g., thickness ratio) are varied. Section \ref{Sec:S4_Shapley_NR_BS} provides the material properties of the rubber. For the third layer, epoxy is used. The Young's modulus and density for epoxy are 3 GPa and 1200 $\mathrm{kg/m^3}$, respectively. Also, the dimensionless thickness of the epoxy layer is taken to be unity. 

The properties of the first layer, serving as the inner core of the sonic crystal, are altered to observe their impact on the design. The parameters subjected to variation include Young's modulus ratio, density ratio, and thickness ratio. These ratios are expressed as $x_{\text{layer1}}:x_{\text{layer2}}$, where $x$ denotes the respective characteristic. It is important to note that both the material and geometrical properties of the third layer remain constant. Consequently, these values are perturbed during the Shapley value analysis.

The QoIs---the objectives of the design---remain the same as for Bragg scattering: achieving a lower first frequency cut-off point and a larger first band-gap. In executing the Shapley value analysis, the base values for all the parameters are taken to unity. Therefore, for these base values, there is no generation of band-gaps, as there are no contrasts in material properties.

\subsection{Data set}
The data were acquired through an eigen-frequency analysis conducted using COMSOL, a finite element software renowned for its multiphysics capabilities \citep{COMSOL}. Following the acquisition of the corresponding QoIs for the parameters, a Python code was employed for Shapley value analysis. This Python code also facilitated data manipulation, specifically tailored for either decreasing the first frequency cut-off point or increasing the width of the first band-gap, separately. The parameters under examination were constrained within the specified ranges outlined in \textbf{Table \ref{Table:dataset_local_range}}. We now discuss the visualization of these data sets. Note that all graphs concerning the elasticity modulus ratio axis were plotted on a logarithmic scale.

\begin{table}[h]
\centering
\caption{Parameter ranges for Shapley value analysis of \emph{sonic crystals}. \label{Table:dataset_local_range}}
\begin{tabular}{l|r}
\centering
\textbf{Parameter} & \textbf{Range}\\
\hline
Young's modulus ratio & 1--100000\\
\hline
Density ratio &  1--10\\
\hline
Thickness ratio & 1--10\\
\end{tabular}
\end{table}

\subsubsection{Visualizing the labeled data for decreasing the first frequency cut-off point.}
\begin{figure}[htp]
    \includegraphics[width=0.45\textwidth]{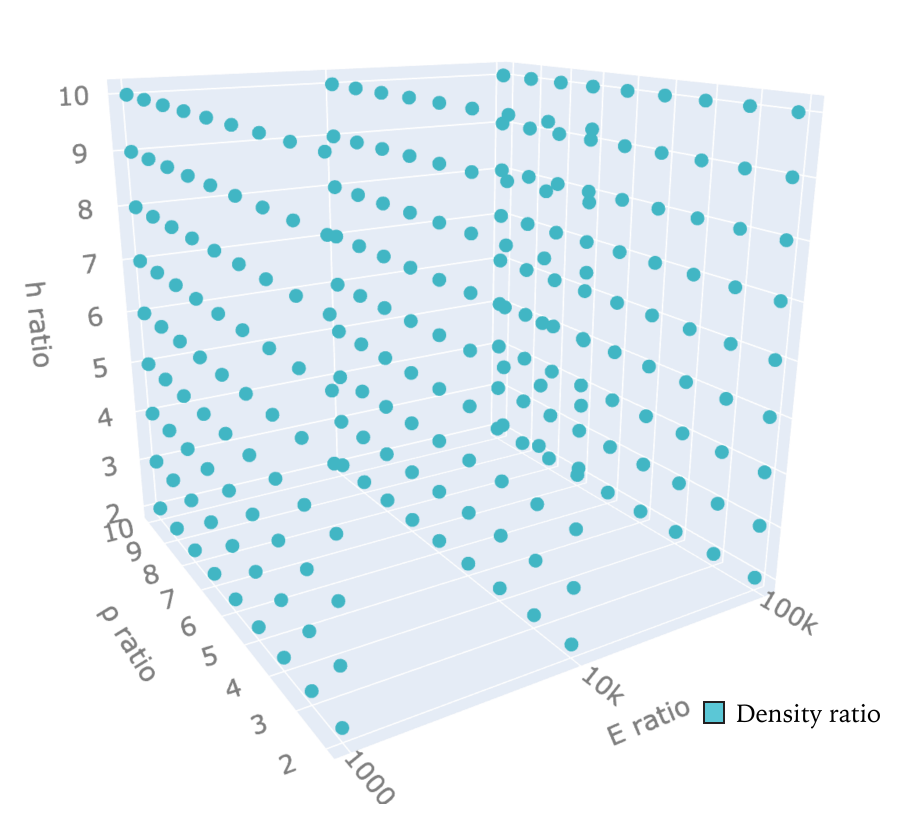}
    \caption{\textsf{Shapley value analysis for a sonic metamaterial.} This figure presents the parameter ranking for the \emph{first} QoI---decreasing the first frequency cut-off point. Over a wide range for the input parameters, the dominant parameter is the density ratio. Other two parameters are not found to be effective. \label{fig:dominantsecondffcoplr}}
\end{figure}

Shapley value analysis was performed to decrease the first frequency cut-off point. As shown in \textbf{Fig. \ref{fig:dominantsecondffcoplr}}, increasing the density ratio---by increasing the density of layer 1---is the only effective parameter. Conversely, increasing Young's modulus and thickness ratios will have the reverse effect, leading to an increase in the first frequency cut-off point.

\subsubsection{Visualizing the labeled data for increasing the first band-gap width.}
    
\begin{figure}[htp]
    \includegraphics[width=1.0\textwidth]{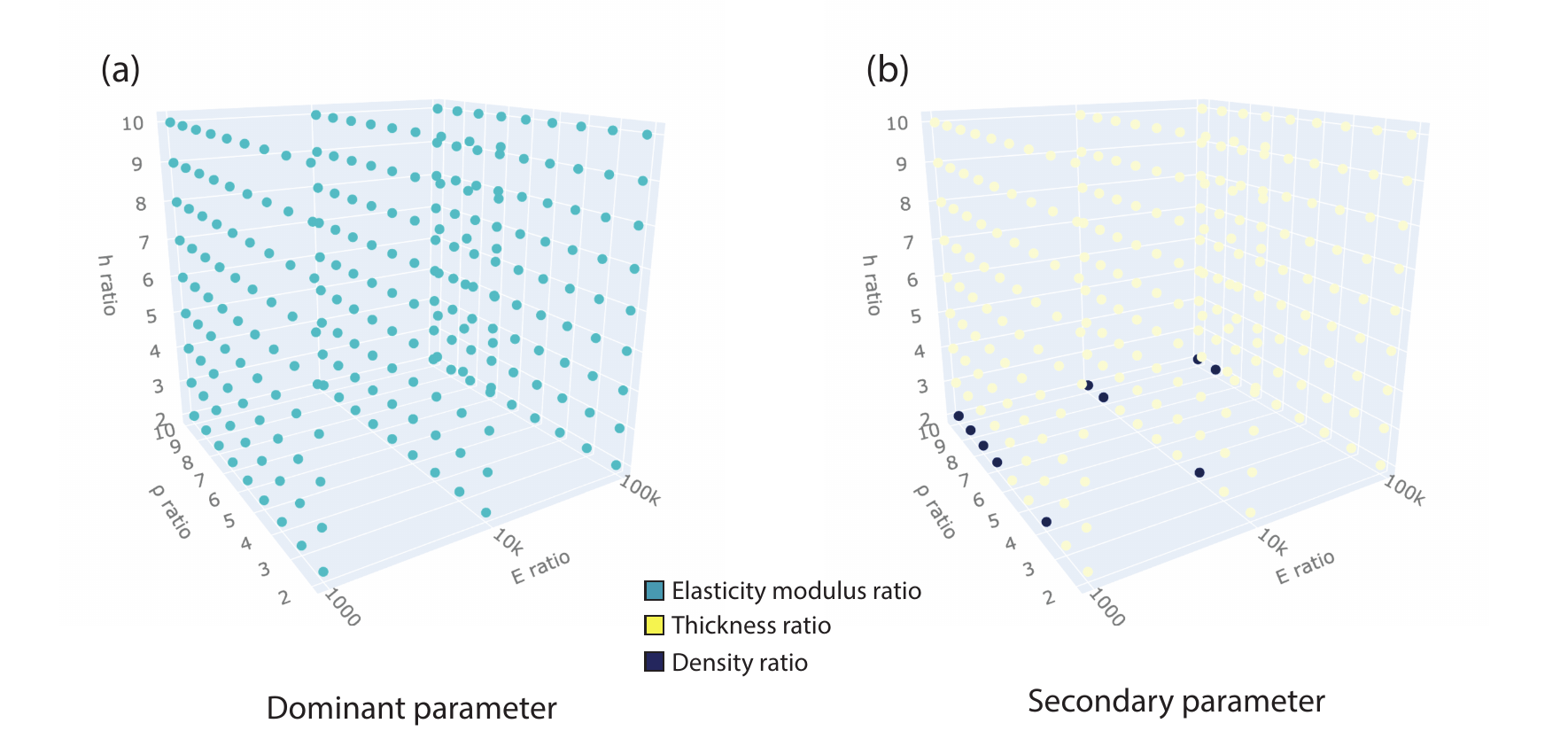}
    \caption{\textsf{Shapley value analysis for a sonic metamaterial.} This figure shows the ranking of effective parameters for the \emph{second} QoI---the first band-gap width. As shown in (a), the dominant parameter is the  Young's modulus ratio for all ranges. For secondary parameter, given in (b), thickness ratio governs the design for most of the ranges. This dominance is replaced by density ratio, when density ratio is higher [$\geq 7$] and thickness ratio is equal to 2.} 
    \label{fig:dominantsecondfbglr}
\end{figure}

Similar to the previous case, Shapley value analysis was performed to identify modifications that increase the first band-gap width. As depicted in \textbf{Fig.~\ref{fig:dominantsecondfbglr}}, increasing the Young's modulus ratio is the dominant parameter for achieving this QoI across all ranges.

\subsection{Discussion on EDA for local resonance}
Since the analytical dispersion relationship is unknown for sonic crystals, COMSOL was used to obtain the band structure for various design parameters. Once the dataset was complete, Shapley value analysis was applied to determine the dominant design parameter. The main observations from visualization and Shapley value analysis are as follows.

To \emph{decrease the first frequency cut-off point},  increasing the density ratio is the only effective approach. Since increasing Young's modulus and thickness ratios tend to increase the first frequency cut-off point, these two ratios should be kept at a minimum as much as permitted by design. 

To \emph{increase the width of band-gap}, Young's modulus is the dominant design parameter for all ranges, followed by thickness and density ratios for different ranges. Although changing any of these three parameters would result in a larger band-gap, the most effective approach is to make the inner core stiffer.

\section{MACHINE LEARNING AND DEEP LEARNING MODELS}
\label{Sec:S5_Shapley_ML}

Supervised learning is a training approach that enables a model to identify an output class---from a predetermined classes of outputs---for newly introduced data. Regression and classification problems often use this approach, as in these problems, the possible output classes are known \emph{a priori}, and the model---after training---has to pick an output for a new data point. However, in clustering problems, the output classes are not known in advance---hence, these problems avail the unsupervised learning approach to identify similarity among data points and divide the data points into groups or clusters \citep{cunningham2008supervised}. Accordingly, all the machine learning algorithms considered in this paper use supervised learning, as the output classes are known.  

Several machine- and deep-learning algorithms are available in the literature to train models via supervised learning. Machine learning models consist of one layer comprising multiple neurons, while deep learning models consist of several to numerous layers with many neurons in each layer. Because of their inherent ability to capture the underlying complexity of data, neural networks can provide accurate and fast results and handle large complex data sets \citep{bonaccorso2017machine}. 

This paper explores several models to predict QoIs for phononic and sonic crystals. For the regression models, the used modeling approaches are (a) neural network, (b) linear regression, and (c) random forest regressor. One of the reasons for choosing the neural network is to compare deep learning models with machine learning models. Linear regression is a simple and common approach for approximating output values by introducing a function that represents the behavior of the input data \citep{vasilev2019python}. This function could be linear, quadratic, or even higher order. Finally, the random forest regression algorithm constructs $n$ decision trees with bootstrap samples and takes the average of the results to predict the new data \citep{segal2004machine}.

We build the mentioned models for Bragg scattering and local resonance data sets in the following subsections. We evaluate the models using the root-mean-squared error (RMSE) and $R^2$ scores. RMSE measures how close the predicted values are to the ground truth and is a popular evaluation method for regression models \citep{chai2014root}:
\begin{align}
\label{Eqn:RMSE}
    \mathrm{RMSE}(y, \hat{y}) = \sqrt{\frac{1}{n}\sum_{i=0}^{n-1}{({y_i -\hat{y}_i})^2}}
\end{align} 
in which $y$ is the ground truth, $\hat{y}$ is the predicted value, and $n$ is the number of samples. $R^2$ score---also called the coefficient of determination---is another common and effective approach for evaluating regression models \citep{renaud2010robust}: 
\begin{align}
\label{Eqn:R2}
    R^2(y, \hat{y}) = 1-\frac{\sum_{i=0}^{n-1}{({y_i -\hat{y}_i})^2}}{\sum_{i=0}^{n-1}{({y_i -\bar{y}_i})^2}}
\end{align}
in which $\bar{y}$ is the mean of ground truth. An $R^2$ score of 0 means that the model is the same as the mean model, while a $R^2$ score closer to 1 indicates a better fit. 

\subsection{Bragg Scattering}
\textbf{Table} \ref{Table:Dataset} provides details on the data used for regression models. The input comprises parameter ratios, while the output is the QoIs, obtained from the dispersion relation. This data set consisted of 2269 data points. 

\begin{table}[h]
\centering
\caption{Data set for Bragg scattering models.} \label{Table:Dataset}
\begin{tabular}{ccccc}
\multicolumn{3}{c|}{Parameter}                                        & \multicolumn{2}{c}{Quantity of Interest} \\ \hline
\multicolumn{1}{c|}{\textbf{E ratio}} & \multicolumn{1}{c|}{\textbf{\begin{math}\rho\end{math} ratio}} & \multicolumn{1}{c|}{\textbf{h ratio}} & \multicolumn{1}{c|}{\thead{First frequency \\ cut-off point}} & \thead{First band-gap} \\ \hline
\multicolumn{1}{c|}{0.1--50000} & \multicolumn{1}{c|}{0.1--9.5} & \multicolumn{1}{c|}{0.1--11} & \multicolumn{1}{c|}{\makecell{From \\ dispersion \\relation}} &  \makecell{From \\ dispersion \\relation}  \\
\end{tabular}
\end{table}

For the neural networks, the input values within the data set---parameter ratios---were scaled using Scikit-learn \citep{sklearn_api}. There is no need to normalize the data for random forest regressor, as it is a tree-based algorithm. This aspect was evident even in our numerical simulations---the results were the same with and without normalized data for the random forest regressor. We also observed that the linear regression gave the same results with and without normalization of the data, although linear regression is not graph-based. For all the models, the data set is divided into train/validation/test subsets. For the reproducibility of the results, the randomness was fixed for all the models. The same seed number was used for the random forest regressor and linear regression when dividing the data sets. For the neural network, the optimum seed number was used.

\subsubsection{Regression Models}
\label{Subsubsec_6.1.4}

As data generation by building a model is not time-efficient for metamaterials that lack the dispersion relationship, a machine learning model that can predict QoIs is needed. 

\begin{figure}[htp]
    \includegraphics[width=1.0\textwidth]{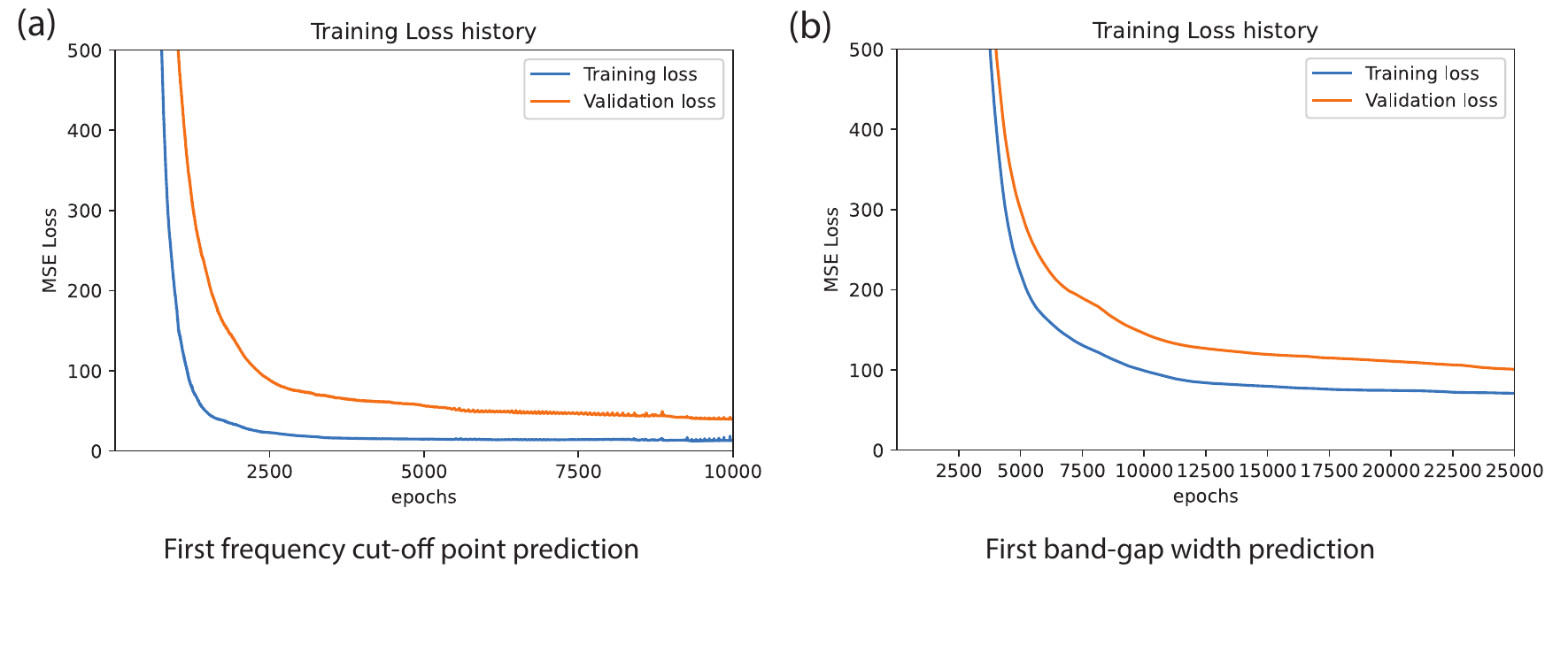}
    \caption{\textsf{Phononic metamaterial.} This figure shows the training and validation loss history. The loss is given in terms of MSE loss. Both graphs indicate over-fitting. Due to the large MSE losses for lower epochs, both graphs are plotted for lower ranges of the loss.}
    \label{fig:loss_hist_bs_regression}
\end{figure}

Two \emph{neural networks}, each comprising two hidden layers, were constructed for regression tasks, employing ReLU as the activation function. The first neural network predicted the initial frequency cut-off point, while the second predicted the width of the first band-gap. This separation facilitates hyper-parameter tuning and accuracy comparison with other models. Adam optimizer was chosen for optimization, while mean squared error loss was utilized as the loss function for these neural networks. For both neural network models, the entire dataset was divided into 80\% for training and 20\% for testing. The validation set was then created by further partitioning the training set into 80\% for training and 20\% for validation. For the model predicting the \emph{first frequency cut-off point}, a learning rate of 0.0025 was selected, with the number of epochs set to 10,000.
Additionally, a weight decay of 0.3 was incorporated. \textbf{Figure~\ref{fig:loss_hist_bs_regression}\textcolor{blue}{a}} depicts The loss histories for both the training and validation subsets. The graph of the neural network's loss history suggests overfitting, with losses presented in terms of the mean squared error (MSE). Subsequently, the final testing root mean squared error (RMSE) value was calculated as 10.49, accompanied by an $R^2$ score of 0.9971. For the second neural network model, predicting the \emph{width of the band-gap prediction}, the learning rate was set to  0.0008, while the number of epochs to 25000. For this model, weight decay was set to 0.5. \textbf{Figure~\ref{fig:loss_hist_bs_regression}\textcolor{blue}{b}}, shows the loss history of training and validation sets and indicates over-fitting. The RMSE and $R^2$ scores for the testing subset were obtained as 14.93 and 0.9993, respectively.

For the \emph{linear regression} models, two separate models were constructed to predict QoIs. These models were trained by incorporating polynomial features. The \emph{first frequency cut-off point prediction} model achieved its optimal performance with third-degree polynomials. The training and validation losses were 81.21 and 81.66, respectively, measured in terms of RMSE loss. The RMSE value for the testing subset was 76.45, with an $R^2$ score of 0.8318, indicating a lack of proper fit to the data. For the \emph{width of band-gap prediction}, a polynomial degree of four was selected. However, the training loss was 131.49, and the validation loss was 135.42, indicating overfitting. The RMSE for the testing subset was 129.42, with an $R^2$ score of 0.9550.

\begin{figure}[htp]
    \includegraphics[width=1.0\textwidth]{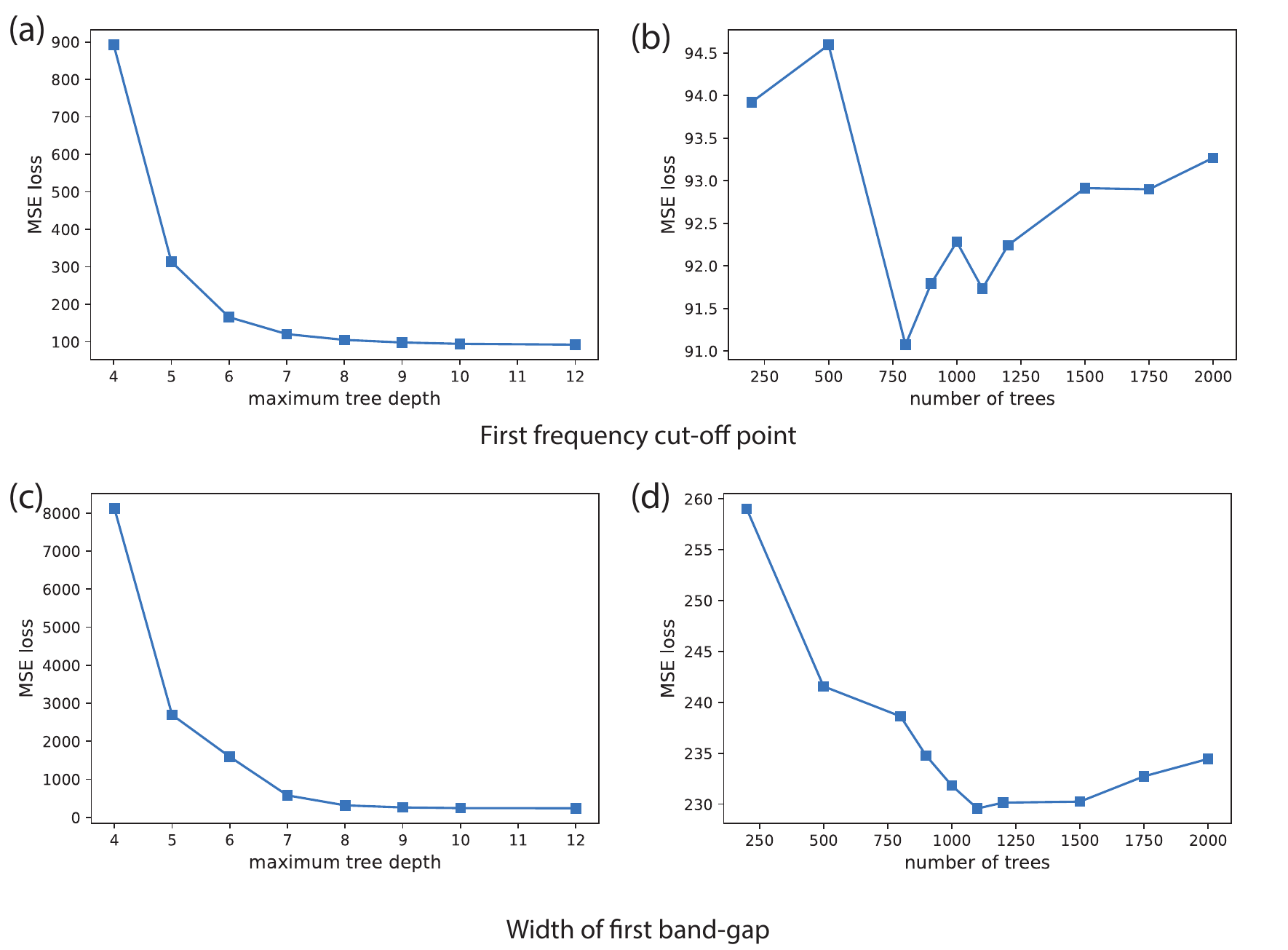}
    \caption{\textsf{Phononic metamaterial.} This figure shows the hyper-parameter tuning of the  \emph{random forest regressor algorithm}. For the first frequency cut-off point,  the optimum maximum tree depth is 10, whereas the optimum number of trees is 800; see subfigures (a) and (b). For the first band-gap width, the optimal maximum tree depth and number of trees are 10 and 1100, respectively; see subfigures (c) and (d).
    \label{fig:loss_hist_bs_hyper}}
\end{figure}

For both the \emph{random forest regressor} models, hyperparameter tuning was conducted through cross-validation. Each model, targeting different QoIs, was individually constructed and fine-tuned for optimal tree depth and the number of estimators. In the case of the \emph{first frequency cut-off point prediction} model, the variation in the loss was plotted against the tuned hyperparameters, depicted in \textbf{Fig.~\ref{fig:loss_hist_bs_hyper}\textcolor{blue}{a-b}}. Achieving a testing RMSE value of 4.16 alongside an impressive $R^2$ score of 0.9996 was possible with the optimum model, characterized by a maximum tree depth of 10 and 800 selected estimators. Similarly, for the \emph{width of band-gap prediction model}, hyperparameter tuning was executed with identical parameters. The tuning process is illustrated in Fig.~\ref{fig:loss_hist_bs_hyper}\textcolor{blue}{c-d}. Subsequent testing revealed an RMSE of 12.75 and an $R^2$ score of 0.9996. In this scenario, the optimal model also featured a maximum tree depth of 10, with the number of estimators set to 1100.

\subsection{Local resonance} The data set for local resonance was constructed using COMSOL. \textbf{Table} \ref{Table:lr_dataset} provides details on this data set. Note that Young's modulus was plotted on a logarithmic scale. In the regression dataset, the inputs comprise parameter values, while the outputs correspond to the QoIs. 

\begin{table}[h]
\centering
\caption{Data set for local resonance models.}
\label{Table:lr_dataset}
\begin{tabular}{c|c|c|l|l}
\textbf{E ratio} &\textbf{\begin{math}\rho\end{math} ratio} & \textbf{h ratio} & \thead{First frequency \\ cut-off point} & \thead{First band-gap}\\
\hline
1--100000 & 1--10 & 1--10  & \makecell{Obtained from \\ COMSOL model} & \makecell{Obtained from \\ COMSOL model} \\
\end{tabular}
\end{table}

For the neural networks, only the input values from the data set---the parameter ratios but not the outputs (QoI values)---were scaled with the same procedure as the one used for Bragg scattering. As indicated for the Bragg scattering models, the random forest regressor and linear regression models used the unnormalized data set. The data set was divided into train/validation/test subsets for all the models. The randomness in the models for the local resonance models was fixed following the procedure used for models under the Bragg scattering.

\subsubsection{Regression Models}
Given that the dispersion relation remains elusive for sonic metamaterials, the construction of computational models becomes imperative for its determination. However, this endeavor proves time-consuming, demanding numerous realizations to attain reliable results. To streamline this process, there arises a pressing need for a faster and more accurate model, a need that machine learning approaches can effectively address. Various regression models have been developed for this purpose, each offering distinct efficiencies, elucidated through their respective RMSE and $R^2$ scores.

\begin{figure}[h]
    \includegraphics[width=1.0\textwidth]{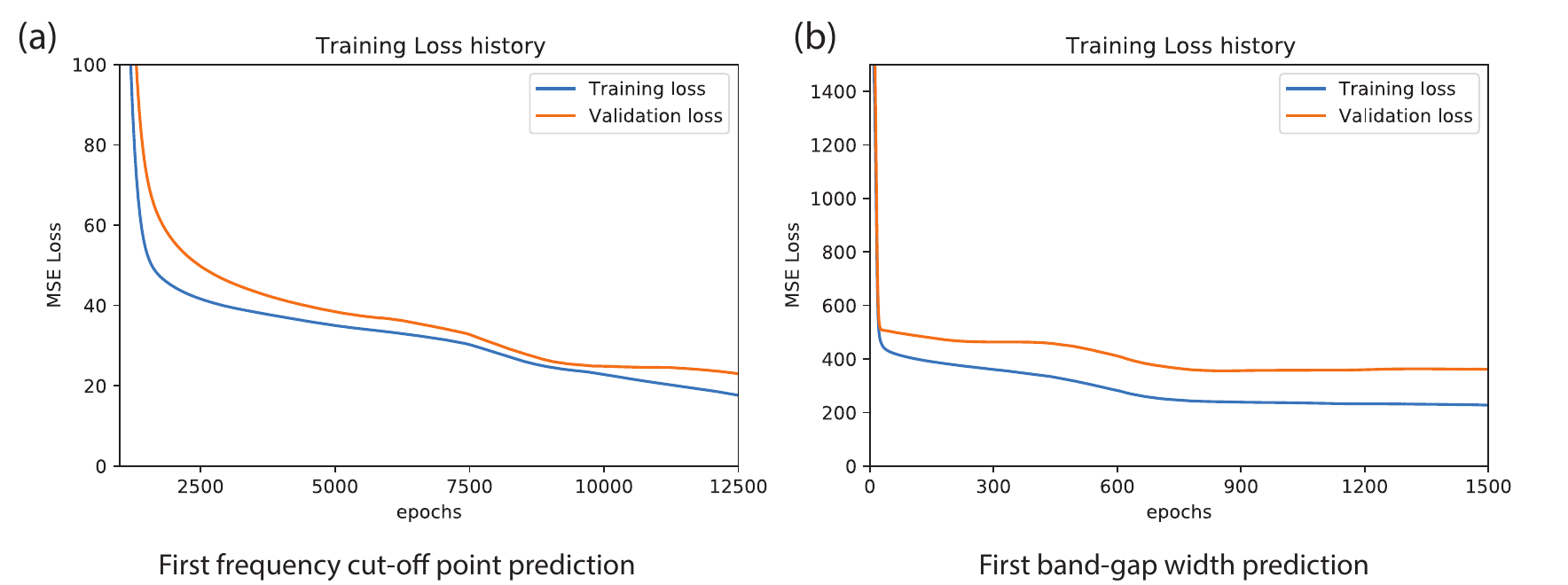}
    \caption{\textsf{Sonic metamaterial.} This figure shows the loss history of training and validation subsets for neural network models. In these graphs, the loss is plotted in terms of MSE loss. The graph for the width of the first band-gap (right) indicates over-fitting.}
    \label{fig:loss_hist_lr_regression}
\end{figure}

Two \emph{neural networks} were constructed separately to predict QoIs. Both models utilized ReLU as the activation function and employed MSE loss as the selected loss function, with the final loss presented as the square root of MSE loss. This loss metric was computed for each subset—training, validation, and testing. Both models were trained and tested with a split of 0.8 for training and 0.2 for testing, while the validation set was derived from the training data with the same split ratio. Each neural network comprised a two-hidden-layer structure. For the \emph{first frequency cut-off point prediction}, the number of epochs was set to 12500. Utilizing the Adam optimizer, the optimal learning rate of 0.0002 was determined, coupled with the utilization of the entire training batch during model training. The training and validation loss history is depicted in \textbf{Fig.~\ref{fig:loss_hist_lr_regression}\textcolor{blue}{a}}. Subsequently, the testing RMSE and $R^2$ score were computed as 4.57 and 0.8395, respectively. On the other hand, for the \emph{width of the band-gap prediction}, the number of epochs was set to 1500. Stochastic gradient descent (SGD) was employed as the optimizer with a learning rate of 0.0002. Mini-batches consisting of 64 data points were used during training. The loss history of both the training and validation subsets is illustrated in Fig.~\ref{fig:loss_hist_lr_regression}\textcolor{blue}{b}. Notably, this graph suggests overfitting. The testing RMSE was determined to be 17.11, accompanied by an $R^2$ score of 0.3663, indicating a poor fit for this model.

Separate \emph{linear regression} models were constructed for each QoI. Optimal performance was achieved by augmenting both models with third-degree polynomial features. Subsequently, the square root of the mean square error (RMSE) was computed for each subset. For the \emph{first frequency cut-off point prediction} model, the training RMSE loss was 3.71, with a validation RMSE of 3.21. On the testing subset, the RMSE value was 3.76, accompanied by an $R^2$ score of 0.8996. Further, for the \emph{band-gap width prediction} model, the training RMSE loss was 9.57, with a validation RMSE of 10.2. The testing RMSE was 8.57, with an $R^2$ score of 0.8893. Notably, unlike the Bragg scattering model, these models exhibited neither overfitting nor compromised accuracy.

\begin{figure}[h]
    \includegraphics[width=1.0\textwidth]{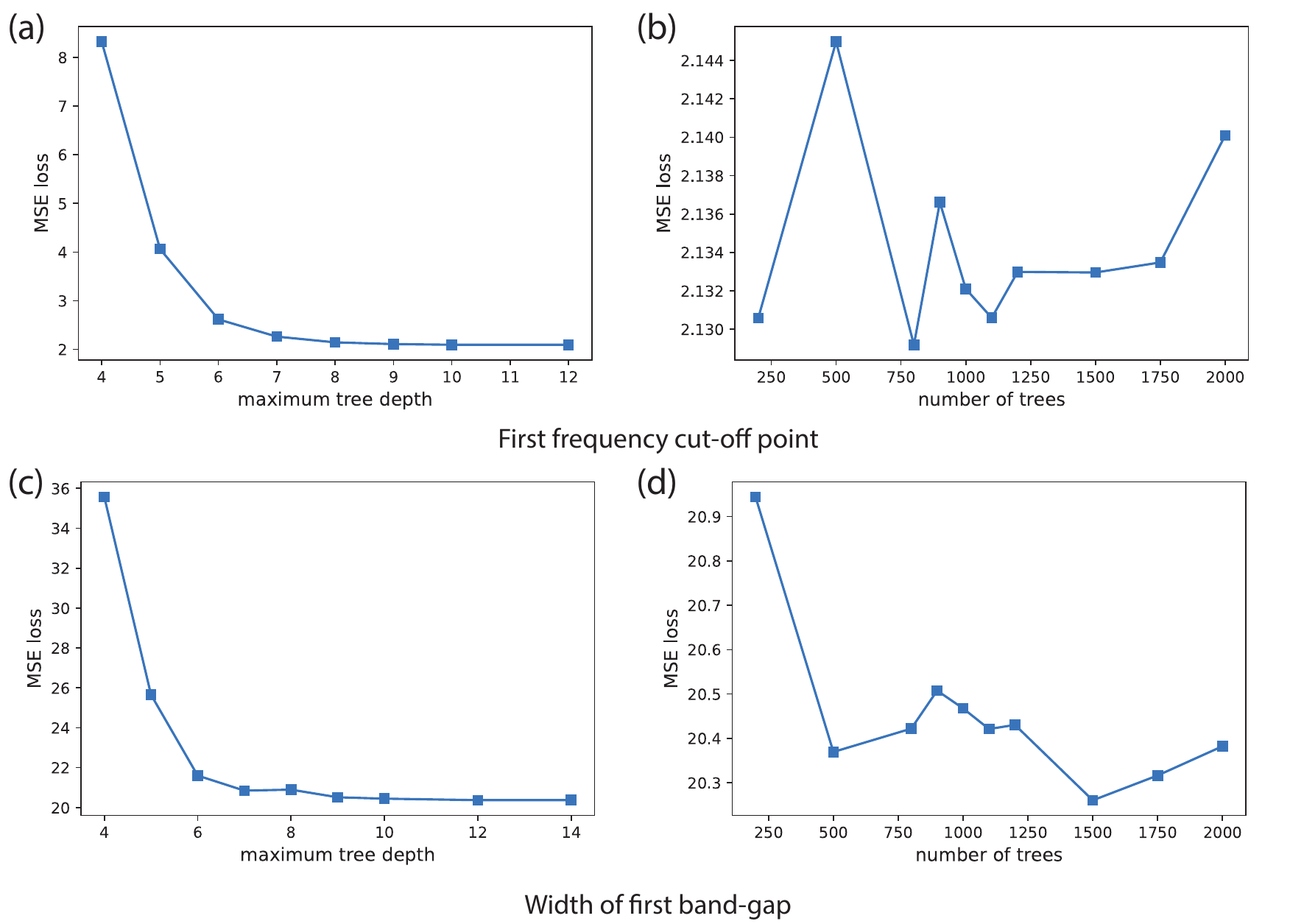}
    \caption{\textsf{Sonic metamaterial.} This figure shows the hyper-parameter tuning of the random forest regressor model. For the first frequency cut-off point, the optimum maximum tree depth is 8, whereas the optimum number of trees is 800; see subfigures (a) and (b). However, the number of trees has no significant effect on the mean squared error (MSE) loss. For the width of the first band-gap, the optimal maximum tree depth and number of trees are 12 and 1500, respectively; see subfigures (c) and (d). 
    \label{fig:loss_hist_lr_hyper}}
\end{figure}

Two \emph{random forest regressor} models were developed using cross-validation and hyper-parameter tuning, focusing on optimizing the maximum tree depth and the number of estimators (i.e., number of trees). The results of the hyper-parameter tuning, illustrating the changes in the loss for both models, are depicted in \textbf{Fig.~\ref{fig:loss_hist_lr_hyper}}, with loss values presented as RMSE. In the \emph{first frequency cut-off point prediction}, a model configured with a maximum tree depth of 8 and 800 estimators achieved an RMSE value of 1.65 and an $R^2$ score of 0.9825 on the test set. On the other hand, in the \emph{band-gap width prediction}, employing a model with a maximum tree depth of 12 and 1500 estimators resulted in an RMSE of 3.55, with an associated $R^2$ score of 0.9827.

\subsection{Discussion on machine learning models}
\textbf{Table \ref{Tab:RMSE_R2}} summaries the obtained RMSE values and $R^2$ scores.

The main findings of the machine learning algorithms for \emph{layered phononic} metamaterials are as follows: For predicting the \emph{first frequency cut-off point}, both neural network and random forest regressor algorithms perform well on the dataset. However, random forest yields a lower RMSE value and is computationally faster than the neural network. On the other hand, linear regression does not provide a reliable prediction for this dataset. For predicting the \emph{band-gap width}, both random forest regressor and neural network models demonstrate good accuracy. However, the random forest model provides the best fit with the least prediction time. In contrast, linear regression exhibits the poorest fit.

The main findings of the machine learning algorithms for \emph{sonic} metamaterials are as follows: For predicting the \emph{first frequency cut-off point}, all three models yield small losses. However, based on their $R^2$ scores, the random forest regressor is the best model. On the other hand, both linear regression and neural network models perform similarly in fitting the data, with linear regression resulting in a slightly better fit. For predicting the \emph{band-gap width}, the best-performing model is the random forest regressor, followed by linear regression. The $R^2$ score and RMSE value indicate that the neural network provides a poor fit.

\begin{table}[h]
\centering
\addtolength{\tabcolsep}{-0.5pt}
\caption{Testing root-mean-square error (RMSE) values and $R^2$ scores of the regression models.}
\label{Tab:RMSE_R2}
\begin{tabular}{c|c|c|c|c|c}
\centering
\textbf{\makecell{Mechanism}} & \textbf{\makecell{Quantity of\\ Interest}} & \textbf{\makecell{Evaluation \\ metric}}  & \textbf{\makecell{Neural \\network}} & \textbf{\makecell{Linear \\regression}} & \textbf{\makecell{Random Forest\\ Regressor}}\\
\hline
\multirow{4}{*}{\textbf{\makecell{Bragg \\scattering}}} & \multirow{2}{*}{\textbf{\makecell{First frequency \\cut-off point}}} & \makecell{RMSE} & 10.49 & 76.45 & 4.19 \\
\cline{3-6}
 & & \makecell{$R^2$ score} & 0.9971 & 0.8318 & 0.9996 \\
\cline{2-6}
& \multirow{2}{*}{\textbf{\makecell{First band-gap \\width}}} & \makecell{RMSE} & 14.93 & 129.42 & 12.75 \\
\cline{3-6}
& & \makecell{$R^2$ score} & 0.9993 & 0.9550 & 0.9996 \\
\hline
\multirow{4}{*}{\textbf{\makecell{Local \\resonance}}} & \multirow{2}{*}{\textbf{\makecell{First frequency \\cut-off point}}} & \makecell{RMSE} & 4.57 & 3.76 & 1.65 \\
\cline{3-6}
 & & \makecell{$R^2$ score} & 0.8395 & 0.8996 & 0.9825 \\
\cline{2-6}
& \multirow{2}{*}{\textbf{\makecell{First band-gap \\width}}} & \makecell{RMSE} & 17.11 & 8.57 & 3.55 \\
\cline{3-6}
& & \makecell{$R^2$ score} & 0.3663 & 0.8893 & 0.9827 \\
\end{tabular}
\end{table}

\section{SUMMARY AND CONCLUDING REMARKS}
\label{Sec:S6_Shapley_CR}
This paper addressed the influence of material and geometrical properties of a metamaterial on band structure characteristics, particularly under the phenomena of Bragg scattering and local resonance. Given our focus on seismic wave barriers, we considered the frequency at which the first cut-off occurs and the magnitude of the first band-gap as the two primary quantities of interest (QoIs). Our contributions are twofold. \emph{Firstly}, we presented a sensitivity analysis method utilizing Shapley values---a technique from cooperative game theory---to quantify the dominance of input parameters on these quantities of interest. \emph{Secondly}, we introduced machine learning tools for obtaining regression models to predict these quantities of interest.

\subsection{Shapley value analysis} 
The main advantages of Shapley value analysis are: (i) it requires a lesser amount of data to rank the input parameters compared to Sobol analysis, (ii) it demarcates the transition regions between dominant parameters, and (iii) it provides information on how to alter the parameters, whether to increase or decrease them.
On the other hand, the main disadvantages are: (i) when compared to Sobol analysis, Shapley value analysis does not rank the interactions among the individual parameters, and (ii) it does not offer regression expressions---or the so-called reduced-order models---for predicting the QoIs.

Based on a Shapley value analysis, the main findings concerning \emph{phononic} metamaterials are as follows:
\begin{enumerate}
\item \textbf{Decreasing the first frequency cut-off point:} Young's modulus ratio does not significantly contribute to this reduction; instead, increasing the ratio tends to raise the first frequency cut-off point. Regarding the other two parameters—thickness and density ratios—within the ranges of [0.1–2] for the density ratio and [0.1–9.5] for the thickness ratio, the thickness ratio demonstrates dominance. Beyond these ranges, the density ratio governs the design and extends its dominance range.
\item \textbf{Increasing the band-gap width:} For small ranges of the thickness ratio ([0.1–1.5]), neither Young's modulus nor density ratios exhibit significant effectiveness. In other words, marginal increases in these parameters tend to reduce the width of the first band-gap. To achieve effectiveness, a larger increase in these parameters is required. For higher values of the thickness ratio, Young's modulus begins to exert similar dominance to that of the thickness ratio. In cases where a substantial increase in Young's modulus is impractical, maintaining this material property while introducing only geometrical changes remains a viable option.
\end{enumerate}

The corresponding conclusions for \emph{sonic} metamaterials are as follows:
\begin{enumerate}
\item \textbf{Decreasing the first frequency cut-off point:} The only effective parameter for achieving this is the density ratio, where a denser inner core is advantageous. Increasing the ratios of the other two parameters---Young's modulus and thickness---has a detrimental effect on this QoI.
\item \textbf{Increasing the width of the first band-gap:} Contrast in Young's modulus, with the inner core being stiffer, emerges as the dominant parameter across all ranges. However, unlike in phononic metamaterials, increasing the thickness ratio ($h_{\text{layer1}}:h_{\text{layer2}}$) or density ratio ($\rho_{\text{layer1}}:\rho_{\text{layer2}}$) also augments this QoI, although these effects are secondary compared to the impact of Young's modulus ratio.
\end{enumerate}

These findings are valuable as they offer insights into how design parameters can be economically adjusted. For instance, when two materials are readily available, a design engineer can manipulate the geometrical properties within specific ranges to influence the design, thereby leveraging existing materials effectively. This approach provides a practical method to manipulate QoIs using available materials.

\subsection{Machine learning}
Three regression models based on neural networks, linear regression, and random forest regressor algorithms were constructed for each QoI, separately analyzing Bragg scattering and local resonance. \emph{For Bragg scattering}, the random forest regressor emerged as the optimal model for both QoIs, while linear regression exhibited the least favorable fit. \emph{Concerning local resonance}, superior evaluation metrics were achieved with the random forest regressor, although all models demonstrated low RMSE values. Notably, the neural network models produced the lowest $R^2$ scores for both QoIs. The chosen machine learning algorithms, except linear regression, exhibited similar performance across the datasets for both phononic and sonic crystals. Given the rapid evolution of the deep learning field, leveraging newer algorithms and computer architectures holds the potential to further reduce computational time. In times to come, these deep-learning-based models could emerge as cost-effective alternatives to conventional numerical formulations.

A plausible future direction could involve developing a framework using deep learning tools to tackle the inverse problem---also known as band engineering. Said differently, what is needed is a deep learning framework for designing metamaterials, which involves selecting base materials and defining geometrical properties to achieve a desired band structure, such as the first frequency cut-off and first band-gap.

\bibliographystyle{plainnat}
\bibliography{Master_References}
\end{document}